\documentclass[11pt, reqno]{amsart}
\usepackage{amssymb, color}
\numberwithin{equation}{section}

\begin{document}
\setlength{\oddsidemargin}{0cm}
\setlength{\evensidemargin}{0cm}
\baselineskip=20pt

\title[$\mathcal{O}$-operators on Leibniz algebras]
{$\mathcal{O}$-operators and related structures on Leibniz
algebras}

\author{Qinxiu Sun}
\address{Department of Mathematics, Zhejiang University of Science and Technology, Hangzhou, Zhejiang 310023, China}
\email{qxsun@126.com}

\author{Naihuan Jing$^*$}
\address{Department of Mathematics,
   North Carolina State University,
   Raleigh, NC 27695, USA}
\email{jing@math.ncsu.edu}  

\thanks{*Corresponding author: Naihuan Jing}

\begin{abstract}
An $\mathcal{O}$-operator has been used to extend a Leibniz algebra by
its representation. In this paper, we investigate several structures
related to $\mathcal{O}$-operators on Leibniz algebras and introduce
(dual) $\mathcal{O}$N-structures on Leibniz algebras associated to their
representations. It is proved that $\mathcal{O}$-operators and dual
 $\mathcal{O}$N-structures generate each other under certain conditions. It
is also shown that a solution of the strong Maurer-Cartan equation
on the twilled Leibniz algebra gives rise to a dual
$\mathcal{O}$N-structure. Finally, $r-n$ structures, RBN-structures
and $\mathcal{B}N$-structures on Leibniz algebras are thoroughly
studied and their interdependent relations are also studied.
\end{abstract}

\maketitle

\section{Introduction}

Leibniz algebras were introduced by Loday \cite{L, LT}
in the study of periodicity in algebraic K-theory by forgetting the anti-symmetry in Lie algebras. Numerous works have
been devoted to various aspects of Leibniz algebras in both mathematics and physics \cite{BW, C, DMS, DW, VKO,
KM, KW} generalizing many results from Lie algebras and associative algebras.

Rota-Baxter operators were first introduced by Baxter in his study
of fluctuation theory in probability \cite{Ba}. It is known that
there is a close relationship between Rota-Baxter operators and the
classical Yang-Baxter equation on Lie algebras, as the former is an
operator form of a classical $r$-matrix \cite{S} under some
conditions. Rota-Baxter operators on associative algebras have been
found useful in several contexts, for example in quantization of
Poisson geometry, see \cite{G, ZBG} for more information.

The notion of $\mathcal{O}$-operators (also called generalized
Rota-Baxter operators or relative Rota-Baxter operators) was
originally introduced to generalize the well-known classical
Yang-Baxter equation (CYBE) on Lie algebras \cite{K},
$\mathcal{O}$-operators were known to be helpful in providing
solutions of the CYBE on a larger Lie algebra \cite{B1}. Inspired by
Poisson structures, Uchino \cite{U1} introduced the notion of
generalized Rota-Baxter operators (i.e. $\mathcal{O}$-operators in
\cite{LBS} ). The $\mathcal{O}$-operators on associative algebras
give rise to
dendriform algebras, which are important in studying 
bialgebra theory \cite{B2} and operads \cite{BBGN}. Recently, Sheng
and Tang \cite{ST} introduced the notion of
relative Rota-Baxter
operators (called $\mathcal{O}$-operators throughout the paper) on Leibniz algebras
and studied their applications to Leibniz-dendriform algebras and solutions of the classical Leibniz Yang-Baxter equation.
In this picture the Leibniz-dendriform
algebra captures the essential algebraic structure underlying an
$\mathcal{O}$-operator, which in turn leads to a Leibniz algebra structure
on itself and most importantly a solution of the classical Leibniz
Yang-Baxter equation. Moreover, the twilled Leibniz algebras have also been
considered and Maurer-Cartan elements of the associated graded
Lie algebra (gLa) were given therein.

Nijenhuis operators on Lie algebras have been studied in \cite{D}
and \cite{FF}. In the perspective of deformations of Lie algebras,
Nijenhuis operators canonically give rise to trivial deformations
\cite{NR}. Nijenhuis operators have also been studied on pre-Lie
algebras \cite{WSBL} and Poisson-Nijenhuis structures appeared in
completely integrable systems \cite{MM} and were further studied in
\cite{KM, KR}. The $r-n$ structure over a Lie algebra was studied in
 \cite{RAH}, which is the infinitesimal right-invariant integrating the Lie
algebra Poisson-Nijenhuis structure on the Lie group. Recently, Hu,
Liu and Sheng \cite{HLS} studied the (dual) $\mathcal{O}$N-structure
as a generalization of the $r-n$ structure. The notion of an
$\mathcal{O}$N-structure on a bimodule over an associative algebra
was also studied in \cite{LBS} as an associative
analogue of the Poisson-Nijenhuis structure. Moreover, the
compatible algebraic structures underlying $\mathcal{O}$N-structures
were also considered in \cite{LBS}.

In the current work, 
we would like to generalize the algebraic structures related to
$\mathcal{O}$-operators from Lie algebras and associative algebras
to Leibniz algebras, aiming to study their
 relations with
analogous Rota-Baxter operators and formulate criteria in terms of similar structures of the classical Yang-Baxter equations.
 Realizing that
 many results are closely connected in several situations, we would like to carefully formulate them in hoping to
 show similar results in various contexts. 
 Since representation theory of Leibniz algebras is more
 subtle than that of Lie algebras and the classical Leibniz Yang-Baxter equation $[[r,r]]=0$ is
 different from the classical Yang-Baxter equation 
 related to Lie algebras,
 Leibniz algebras should be studied as an independent
 algebraic structure from Lie algebras, not just as a generalization of
 Lie algebras.

The paper is organized as follows. In Section 2, we give some
elementary results on Leibniz algebras. The (dual-)Nijenhuis pair is
studied in Section 3, and it is noted that a Nijenhuis pair can
generate a trivial deformation. In Section 4, we consider (dual)
 $\mathcal{O}$N-structures and explore their properties, and we show that $\mathcal{O}$-operators and dual
 $\mathcal{O}$N-structures can generate each other under some
conditions. In Section 5, we prove that a solution of the strong
Maurer-Cartan equation on the twilled Leibniz algebra affords 
 a dual $\mathcal{O}$N-structure. Finally,
 the notions of $r-n$ structures, RBN-structures and $\mathcal{B}N$-structures on
Leibniz algebras and their relationships are studied in Section 6.

Let $k$ be a field of
characteristic 0 throughout the paper. Unless otherwise specified, all vector spaces are
finite dimensional over $k$.

\section{Preliminaries on Leibniz
algebras }

We briefly review some elementary notions on Leibniz
algebras \cite{L, ST, B, FM}.

A {\it Leibniz algebra} $(\mathfrak g,[ \ , \ ])$ is a vector space
equipped with a bilinear product $[ \ , \ ]: \mathfrak g\times
\mathfrak g\longrightarrow\mathfrak g$ satisfying the (left) Leibniz
property: for any $x_0, x_1, x_2\in\mathfrak g$
\begin{equation*} [x_0, [x_1, x_2]]=[[x_0, x_1], x_2]+[x_1,[x_0, x_2]].
\end{equation*}
Therefore a Lie algebra is certainly a Leibniz algebra.

Let $V$ be a vector space, the endomorphism algebra
$\mathrm{End}(V)$ becomes the Lie algebra $\mathfrak{gl}(V)$ under
the usual bracket operation: $[f, g]=fg-gf$ for $f,
~g\in\mathrm{End}(V)$. Let $\rho^L, \rho^R: \mathfrak
g\longrightarrow \mathfrak{gl}(V)$ be two linear maps satisfying the
following properties:
\begin{align*}
\rho^{L}([x_0,x_1])&=[\rho^{L}(x_0),\rho^{L}(x_1)], \qquad \rho^{R}([x_0,x_1])=[\rho^{L}(x_0),\rho^{R}(x_1)], \\
\rho^{R}(x_1)\rho^{L}(x_0)&=-\rho^{R}(x_1)\rho^{R}(x_0),~~~\forall~x_0,x_1\in
\mathfrak g,
\end{align*}
then the triple $(V, \rho^L, \rho^R)$ is called a {\it
representation} of the Leibniz algebra $(\mathfrak g, [ \ , \ ])$.
In the Lie algebra case, if $(V, \rho)$ is a Lie algebra representation, then
one can take $\rho^L=\rho$ and
$\rho^R=-\rho$ or $\rho^R=0$ to view $V$ as a Leibniz algebra representation. 
In essence, $V$ is a representation of the Leibniz algebra
$\mathfrak g$ if and only if $\mathfrak g\oplus V$ is a
 Leibniz algebra under the bracket
$$[x_0+w_0, x_1+w_1]=[x_0, x_1]+\rho^L(x_0)w_1+\rho^R(x_1)w_0,
~~\forall~~ w_0, w_1\in V,x_0, x_1\in \mathfrak g.
$$
 This Leibniz algebra is called the {\it semidirect product} of $(\mathfrak g, [ \ , \
 ])$ and $(V, \rho^L, \rho^R)$, denoted as $\mathfrak g\ltimes
 V$.

Let $(V, \rho^L, \rho^R)$ be a representation of a Leibniz algebra $\mathfrak g$,
and let $V^*$ be the dual space of the vector space $V$.
The triple $(V^{*} ,(\rho^{L})^{*} ,
-(\rho^{L})^{*}-(\rho^{R})^{*})$ is a representation of $\mathfrak g$ called the {\it dual representation} of $\mathfrak g$ associated
with $V$. 
Here $(\rho^{L})^{*},~(\rho^{R})^{*}$ are defined by (similar to Lie
algebras): $\langle (\rho^{L})^{*}(y)w^{*},w_0\rangle=-\langle
w^{*}$, $\rho^{L}(y)w_0\rangle,~\langle
(\rho^{R})^{*}(y)w^{*}$, $w_0\rangle=-\langle w^{*},
\rho^{R}(y)w_0\rangle$, for all $y\in\mathfrak g,w^*\in V^*,
w_0\in V$    
(see details in \cite{ST}). We remark that there is another version
of dual representations for Leibniz algebras (cf. \cite{BB}).

Let $(V, \rho^L, \rho^R)$ and $(V', {\rho^L}', {\rho^R}')$ be two
representations of a Leibniz algebra $(\mathfrak g,[ \ , \ ])$. A
linear map $f:V\longrightarrow V'$ is called a homomorphism if it
satisfies
$$f \rho^L(x)={\rho^L}'(x) f,~~f \rho^R(x)={\rho^R}'(x) f,~\forall~ x\in  \mathfrak g.$$   

Let $L$ (resp. $R$) be the left (resp. right) multiplication
operator associated to $(\mathfrak g,[ \ , \ ])$, i.e.
$$L(x_0)x_1 = R(x_1)x_0 = [x_0, x_1],~~ \forall ~x_0,x_1\in
\mathfrak g.$$ Then $(\mathfrak g, L, R)$ is a priori a
representation of $\mathfrak g$ called the {\it regular
representation}. Subsequently,
$(\mathfrak{g}^{*},L^{*},-L^{*}-R^{*})$ is a representation of $\mathfrak{g}$, which
is called the {\it dual regular representation}.

{\bf Definition 2.1.} \cite{AG} Let $(\mathfrak g_{1}, [ \ , \ ]_1)$
and $(\mathfrak g_{2}, [ \ , \ ]_2)$ be two Leibniz algebras.
Suppose that there are linear maps
$\rho^{L}_{1},~\rho^{R}_{1}:\mathfrak g_{1} \longrightarrow
\mathfrak{gl}(\mathfrak g_{2})$ and
$\rho^{L}_{2},~\rho^{R}_{2}:\mathfrak g_{2} \longrightarrow
\mathfrak{gl}(\mathfrak g_{1})$ such that $(\mathfrak
g_{2},\rho^{L}_{1},\rho^{R}_{1})$ is a representation of $\mathfrak
g_{1}$ and $(\mathfrak g_{1},\rho^{L}_{2},\rho^{R}_{2})$ is a
representation of $\mathfrak g_{2}$ and they satisfy the following
conditions:
$$\rho^{R}_{1}(x)[a,b]_2-[a,\rho^{R}_{1}(x)b]_2+[b,\rho^{R}_{1}(x)a]_2-\rho^{R}_{1}(\rho^{L}_{2}(b)x)a+\rho^{R}_{1}(\rho^{L}_{2}(a)x)b=0,$$
$$\rho^{L}_{1}(x)[a,b]_2-[\rho^{L}_{1}(x)a,b]_2-[a,\rho^{L}_{1}(x)b]_2-\rho^{L}_{1}(\rho^{R}_{2}(a)x)b-\rho^{R}_{1}(\rho^{R}_{2}(b)x)a=0,$$
$$[\rho^{L}_{1}(x)a,b]_2+\rho^{L}_{1}(\rho^{R}_{2}(a)x)b+[\rho^{R}_{1}(x)a,b]_2+
\rho^{L}_{1}(\rho^{L}_{2}(a)x)b=0,$$
$$\rho^{R}_{2}(x)[x,y]_1-[x,\rho^{R}_{2}(a)y]_1+[y,\rho^{R}_{2}(a)x]_1-\rho^{R}_{2}(\rho^{L}_{1}(y)a)x+\rho^{R}_{2}(\rho^{L}_{1}(x)a)y=0,$$
$$\rho^{L}_{2}(x)[x,y]_1-[\rho^{L}_{2}(a)x,y]_1-[x,\rho^{L}_{2}(a)y]_1-\rho^{L}_{2}(\rho^{R}_{1}(x)a)y-\rho^{R}_{2}(\rho^{R}_{1}(y)a)x=0,$$
$$[\rho^{L}_{2}(a)x,y]_1+\rho^{L}_{2}(\rho^{R}_{1}(x)a)y+[\rho^{R}_{2}(a)x,y]_1+
\rho^{L}_{2}(\rho^{L}_{1}(x)a)y=0,$$ for any $x, y \in \mathfrak
g_{1}, a, b \in \mathfrak g_{2}$. Then there is a Leibniz algebra
 structure on the vector space $\mathfrak g_{1}\oplus \mathfrak g_{2}$ given by
$$[x+a,y+b]=[x,y]_1+\rho^{L}_{2}(a)y+\rho^{R}_{2}(b)x+[a,b]_2+\rho^{L}_{1}(x)b+\rho^{R}_{1}(y)a$$
We denote this Leibniz algebra by $(\mathfrak g_{1}\bowtie \mathfrak
g_{2}, [ \ , \ ])$. The triple $((\mathfrak g_{1}, [ \ , \ ]_1),(\mathfrak
g_{2},[ \ , \
]_2),\rho^{L}_{1},\rho^{R}_{1},\rho^{L}_{2},\rho^{R}_{2})$
satisfying the above conditions is called a {\it matched pair of
Leibniz algebras}. On the other hand, every Leibniz algebra which is
a direct sum of the underlying vector spaces of two subalgebras can
be obtained in this manner. Here a Leibniz {\it subalgebra} is understood in the usual sense.

If a Leibniz algebra $\mathfrak g$ decomposes itself
into a sum of two subspaces: ${\mathfrak g}=\mathfrak g_{1}\oplus
 \mathfrak g_{2}$ such that $\mathfrak g_1, \mathfrak g_{2}$ are subalgebras of $\mathfrak g$, then
the triple $(\mathfrak g, \mathfrak g_{1}, \mathfrak g_{2})$ is
called a {\it twilled Leibniz
 algebra}. By Proposition 3.14 \cite{ST}, we know that there is a one-to-one correspondence between twilled
Leibniz algebras and matched pairs of Leibniz algebras.

Let
 $(V,\rho^{L},\rho^{R})$ be a representation of a Leibniz algebra $(\mathfrak g,[ \ , \
])$. To afford a Leibniz algebra structure on $V$, we need the
notion of an $\mathcal{O}$-operator introduced by Kupershmidt
\cite{K}. A linear map $K:V\longrightarrow \mathfrak g$ is called an
{\it $\mathcal{O}$-operator} (also called a relative Rota-Baxter
operator \cite{ST}) on a Leibniz algebra $(\mathfrak g,[ \ , \ ])$ associated to a
representation $(V,\rho^{L},\rho^{R})$ if for $u, v\in V$
\begin{equation}\label{ks}
[K(u), K(v)]=K(\rho^{L}(K(u))v+\rho^{R}(K(v))u).
\end{equation}
If one writes $u\lhd^{K} v=\rho^{L}(K(u))v~~ \hbox{and}~~u\rhd^{K}
v=\rho^{R}(K(v))u$, then the condition of \eqref{ks} is $[K(u),
K(v)]=K(u\lhd^K v+u\rhd^K v)$. $V$ becomes a Leibniz algebra
under the bracket: $[u, v]^K:=u\rhd^K v+ u\lhd^K v$, which will be
denoted by $(V_{K},[ \ , \ ]^K)$.

The $\mathcal{O}$-operator on a Leibniz algebra $(\mathfrak g,[ \ ,
\ ])$ associated to its regular representation $(\mathfrak g,L,R )$
is called a {\it Rota-Baxter operator}.

The Leibniz algebra $(V_{K}, [ \ , \ ]^{K})$ has a natural
representation $(\mathfrak g,\varrho_{K}^{L},~\varrho_{K}^{R})$
given by
\begin{equation}\label{e:natural}
\varrho_{K}^{L}(v)x=[K(v),x]-K(\rho^{R}(x)v),~~\varrho_{K}^{R}(v)x=[x,K(v)]-K(\rho^{L}(x)v)
\end{equation}
for any $x\in \mathfrak g,v\in V_K$.

With the given $\mathcal{O}$-operator $K$, the space $(\mathfrak
g\oplus V,\{ \ , \ \}_{K})$ becomes a Leibniz algebra with
$$\{x_0+v_0,x_1+v_1\}_{K}=[x_0,x_1]+\varrho_{K}^{L}(v_0)x_1+\varrho_{K}^{R}(v_1)x_0+\rho^{L}(x_0)v_1+\rho^{R}(x_1)v_0+[v_0,v_1]^{K}$$
for any $x_0,x_1\in \mathfrak g,v_0,v_1\in V$.

A permutation $\sigma\in S_n$ is called a {\it $(j,n-j)$-shuffle} if
$\sigma(1)<
 \cdot\cdot\cdot<\sigma(j)$ and $\sigma(j+1)<
 \cdot\cdot\cdot<\sigma(n)$. If $j=0$ or $n$, we assume that $\sigma= id$. The set of all $(j,n-j)$-shuffles
 is denoted by $S_{(j,n-j)}$. The notion of a $(j_1, \cdot\cdot\cdot ,
j_k)$-shuffle and the set $S_{(j_1, \cdot\cdot\cdot , j_k)}$ of $(j_1, \cdot\cdot\cdot , j_k)$-shuffles
are defined similarly.

Let $\mathfrak g$ be a vector space. Put $C^{n}(\mathfrak
g,\mathfrak g)
 =\hbox{Hom}(\otimes^{n}\mathfrak g,\mathfrak g)$
and $C^{*}(\mathfrak g,\mathfrak g)=\oplus_{n\geq 1}C^{n}(\mathfrak
g,\mathfrak g)$. It is known that $C^{*}(\mathfrak g,\mathfrak g)$
is a graded Lie algebra (gLa) with the Balavoine bracket $\{ \ ,\
\}^{B}$:  
$$\{\varphi_1,\varphi_2\}^{B}=\varphi_1\bar{\circ}\varphi_2-(-1)^{mn}\varphi_2\bar{\circ}\varphi_1,$$
where $\varphi_1\bar{\circ}\varphi_2\in C^{m+n+1}(\mathfrak
g,\mathfrak g)$ is given by
$$\varphi_1\bar{\circ}\varphi_2=\sum_{k=1}^{n+1}\varphi_1\circ_{k}\varphi_2,$$
where $\circ_{k}$ is defined by: 
\begin{eqnarray*}
&&\varphi_1\circ_{k}\varphi_2(x_{1},\cdot\cdot\cdot,x_{m+n+1})\\&=&
\sum_{\sigma\in
S_{(k-1,n)}}(-1)^{\sigma}(-1)^{n(k-1)}\varphi_1(x_{\sigma(1)},\cdot\cdot\cdot,x_{\sigma(k-1)},\varphi_2(x_{\sigma(k)},\cdot\cdot\cdot,x_{\sigma(k+n-1)},
x_{k+n}),x_{k+n+1} ,\cdot\cdot\cdot,x_{m+n+1}),
\end{eqnarray*}
for any $\varphi_1\in C^{m+1}(\mathfrak g,\mathfrak g),\varphi_2\in
C^{n+1}(\mathfrak g,\mathfrak g)$.

 Furthermore, if $(\mathfrak g, \mu=[ \ , \ ] )$ is a Leibniz algebra,
then $(C^{*}(\mathfrak g,\mathfrak g), \{ \ ,  \ \}^{B},d)$ is a
differential graded Lie algebra (dgLa), where the differential $d$
is defined by $d(\varphi)=\{\mu,\varphi\}^{B}$. See \cite{B,FM} for
more details.

 Let $\mathfrak g_1, \mathfrak g_2$ be two vector spaces, the lift of the linear map
 $f:\mathfrak g_{i_1}\otimes \mathfrak g_{i_2}\otimes\cdot\cdot\cdot \otimes \mathfrak g_{i_n}\longrightarrow
\mathfrak g_j$ is the map $\hat{f}: (\mathfrak g_1\oplus\mathfrak
g_2)^{\otimes n}\longrightarrow \mathfrak g_1\oplus\mathfrak g_2$
given by
\[\hat{f}((x_1,a_1)\otimes (x_2,a_2)\otimes\cdot\cdot\cdot\otimes
(x_n,a_n))=\left \{
\begin{array}{ll} (f(\alpha_1\otimes \alpha_2\otimes \cdot\cdot\cdot\otimes
\alpha_n),0),~~~~~~~~if~~~j=1,\\
(0,f(\alpha_1\otimes \alpha_2\otimes \cdot\cdot\cdot\otimes
\alpha_n)),~~~~~~~~if ~~~j=2,
\end{array}
\right.\] where $\alpha_k=x_k $ or $a_k$ according to which one $\in
g_{i_k}$. For information on lifts, see e.g. \cite{ST,
U2}.

Let $(\mathfrak g,[ \ , \ ])$ be a Leibniz algebra and $(V, \rho^L,
\rho^R)$ its representation. Consider the graded vector space
$$C^{*}(V,\mathfrak g)=\oplus_{n\geq 1}C^{n}(V,\mathfrak g)=\oplus_{n\geq 1}\hbox{Hom}(\otimes^{n}V,\mathfrak g).$$

{\bf Theorem 2.2.} \cite{ST} With the above notations,
$(C^{*}(V,\mathfrak g), \{\{ \ ,  \ \}\})$ is a graded Lie algebra,
where $\{\{ \ ,  \ \}\}:C^{m}(V,\mathfrak g)\times C^{n}(V,\mathfrak
g)\longrightarrow C^{m+n}(V,\mathfrak g)$ is given by
$$\{\{ \varphi , \psi
\}\}=(-1)^{m}\{\{\mu_{\ltimes},\hat{\varphi}\}^{B},\hat{\psi}\}^{B}$$
for any $\varphi\in C^{m}(V,\mathfrak g)$, where $\psi\in C^{n}(V,\mathfrak
g)$ and $\mu_{\ltimes}$ is the Leibniz bracket of the semidirect
product $\mathfrak g\ltimes V$. Furthermore, its Maurer-Cartan
elements are relative Rota-Baxter operators on the Leibniz algebra
$(\mathfrak g, [ \ , \ ])$ with respect to the representation $(V,
\rho^L, \rho^R)$.

 Now consider the dual regular representation $(\mathfrak g^{*},L^{*},-L^{*}-R^{*})$ of $\mathfrak g$ in the above theorem.
  For $k\geq 1$,
 define $\Psi:\otimes^{k+1} \mathfrak g\longrightarrow
 \hbox{Hom}(\otimes ^{k}  \mathfrak g^{*}, \mathfrak g)$ and
 $ \Upsilon:\hbox{Hom}(\otimes ^{k}  \mathfrak g^{*}, \mathfrak g)\longrightarrow \otimes^{k+1} \mathfrak g$
 respectively by
$$\langle \Psi(P)(\xi_1,\cdot\cdot\cdot,\xi_k),\xi_{k+1}\rangle=\langle P, \xi_1\otimes\cdot\cdot\cdot\otimes \xi_k \otimes\xi_{k+1}\rangle$$
 and
$$\langle \Upsilon(f),\xi_1\otimes\cdot\cdot\cdot\otimes \xi_{k}\otimes\xi_{k+1}\rangle=\langle f(\xi_1,\cdot\cdot\cdot \xi_k ), \xi_{k+1}\rangle$$
for any $P\in \otimes^{k+1} \mathfrak g,~\xi_1,\ldots, \xi_{k+1}\in
\mathfrak g^{*}$ and $f\in \hbox{Hom}(\otimes ^{k} \mathfrak g^{*},
\mathfrak g)$, then there is a graded Lie bracket $[[\ , \ ]]$
 on $\oplus_{k\geq 2}(\otimes^{k}\mathfrak g)$ given by
$$[[P,Q]]=\Upsilon\{\{\Psi(P),\Psi(Q)\}\},~~\forall P\in \otimes^{m+1}\mathfrak g,Q\in \otimes^{n+1}\mathfrak g.$$

Consider the special case of $P,Q\in \otimes^{2}\mathfrak
 g$.
 An
element $r \in \hbox{Sym}^{2}(\mathfrak g)$ is called a classical
{\it Leibniz r-matrix} if it satisfies the classical Leibniz
Yang-Baxter equation (YBE): $$[[r,r]]=0.$$
 It is known that if $\pi$ is a
symmetric solution of the classical Leibniz YBE, then the map
$\pi^{\sharp} : \mathfrak g^{*} \longrightarrow \mathfrak g$ given
by
$$\langle \pi^{\sharp}(\alpha_0),\alpha_1\rangle=\pi(\alpha_0,\alpha_1),~~\forall ~\alpha_0,\alpha_1\in \mathfrak g^{*}$$
is an $\mathcal{O}$-operator on $\mathfrak g$ with respect to the dual representation
$(\mathfrak g^{*},L^{*},-L^{*}-R^{*})$. For more details about the classical Leibniz Yang-Baxter equation, see \cite{ST}.

\section{Nijenhuis pairs and dual Nijenhuis pairs}

\quad Let $(V,\rho^{L},\rho^{R})$ be a representation of a Leibniz
algebra $(\mathfrak g, [ \ , \ ])$. Suppose $\omega : \mathfrak
g\times \mathfrak g\longrightarrow \mathfrak g$ is a bilinear map,
and let $\varpi^{L},\varpi^{R}: \mathfrak g \longrightarrow
\mathfrak{gl}(V)$ be two linear maps. We can use these maps to
deform the Leibniz bracket $[ \ , \ ]$ and the action maps $\rho^L,
\rho^R$ as follows. For real $t\geq 0$ and $x_0, x_1\in\mathfrak g$,
put
 $$[x_0,x_1]_{t}  = [x_0, x_1]
+ t\omega(x_0, x_1),$$
$$ \rho^{R}_{t}(x_0) =  \rho^{R}(x_0)
+ t\varpi^{R}(x_0),~~\rho^{L}_{t}(x_0) =  \rho^{L}(x_0) +
t\varpi^{L}(x_0),~~\forall~x_0,x_1\in \mathfrak g.
$$
The triple $(\omega, \varpi^L, \varpi^R)$ is said to {\it generate
an infinitesimal deformation} of the Leibniz algebra $(\mathfrak g,
[ \ , \ ])$ with the representation $(V,\rho^{L},\rho^{R} )$ if for
each $t$, $(\mathfrak g, [ \ , \ ]_{t} )$ is a Leibniz algebra and
$(V,\rho^{L}_{t},\rho^{R}_{t} )$ is its representation.
Of course, $( \mathfrak g,\rho^{L},\rho^{R} )$ can be
viewed as its ``trivial'' infinitesimal deformation with $( \omega,
\varpi^L, \varpi^R)=(0, 0, 0)$. 

The following concept generalizes that of Lie algebras studied in \cite{HLS}.

 {\bf Definition 3.1.} Two infinitesimal deformations $(\mathfrak g, [ \ , \ ]_t^{},\rho^{L}_{t},\rho^{R}_{t})$,
 $(\mathfrak g, [ \ , \ ]_t^{'},\rho^{L'}_{t},\rho^{R'}_{t})$ of a Leibniz algebra $(\mathfrak g, [ \
, \ ])$ with a representation $(V,\rho^{L},\rho^{R} )$
are {\it equivalent} if there is a homomorphism $(I_{\mathfrak g} +
tN, I_{V} + tS)$ from  $(\mathfrak g, [ \ , \
]_t^{'},\rho^{L'}_{t},\rho^{R'}_{t})$  to
 $(\mathfrak g, [ \ , \ ]_t^{},\rho^{L}_{t},\rho^{R}_{t})$, where $N\in \mathfrak{gl}(\mathfrak g)$ and
 $S\in \mathfrak{gl}(V)$, such that for any
$x_0,x_1\in \mathfrak g$,
\begin{align}\label{deform1}
(I_{\mathfrak g} + tN)([x_0,x_1]_{t}^{'} )& = [(I_{\mathfrak g} +
tN)(x_0),(I_{\mathfrak g} + tN)(x_1)]_{t} , \\ \label{deform2}
(I_{V} + tS)\rho^{L'}_{t}(x_0)& = \rho^{L}_{t}((I_{\mathfrak g} +
tN)(x_0)) (I_{V} + tS),\\ \label{deform3} (I_{V} +
tS)\rho^{R'}_{t}(x_0)& = \rho^{R}_{t}((I_{\mathfrak g} + tN)(x_0))
(I_{V} + tS).
\end{align}
When $(\mathfrak g, [ \ , \ ]_t,\rho^{L}_{t},\rho^{R}_{t})$ is
equivalent to $(\mathfrak g, [ \ , \ ],\rho^{L},\rho^{R})$, the
former is called a {\it trivial infinitesimal deformation} of the
latter.

 By direct calculation, if $([ \ , \ ]_t,\rho^{L}_{t},\rho^{R}_{t})$ is a
trivial infinitesimal deformation, then \eqref{deform1} implies that
$N$ satisfies the following relation:
\begin{equation}\label{Ni}
[N(x_0), N(x_1)]=N([N(x_0), x_1]+[x_0, N(x_1)]-N[x_0,
x_1]),~~\forall ~x_0,x_1\in \mathfrak g,
\end{equation}
which is usually referred to as a {\it Nijenhuis operator}. We also
denote
\begin{equation}\label{e:Nijen}
[x_0, x_1]_{N} =[N(x_0), x_1]+[x_0, N(x_1)]-N[x_0, x_1].
\end{equation}

 Moreover, for every $y \in \mathfrak g$ and $w \in V$, $N, S$ also satisfy
\begin{align}\label{3.1}
\rho^{L}(N(y)) S(w) &=S(\rho^{L}(N(y))
 w)+S(\rho^{L}(y)S(w))-S^{2}(\rho^{L}(y)w),\\ \label{3.2}
\rho^{R}(N(y)) S(w) &=S(\rho^{R}(N(y))
 w)+S(\rho^{R}(y)S(w))-S^{2}(\rho^{R}(y)w).
 \end{align}

Now we are ready to define Nijenhuis pairs and dual
Nijenhuis pairs.

 {\bf Definition 3.2.} A pair $(N, S)$ with $N
\in\mathfrak{gl}(\mathfrak g),S\in \mathfrak{gl}(V)$ is called
 a {\it Nijenhuis pair} on a Leibniz algebra $(\mathfrak g, [ \ , \ ])$ with
 a representation $(V, \rho^{L},\rho^{R})$ if $N$ is a Nijenhuis
operator on $(\mathfrak g, [ \ , \ ])$ and $(N, S)$ satisfies the
identities \eqref{3.1}--\eqref{3.2}. 

{\bf Definition 3.3.} A pair $(N, S)$ with $N
\in\mathfrak{gl}(\mathfrak g),~S\in \mathfrak{gl}(V)$ is called a
{\it dual-Nijenhuis pair} on a Leibniz algebra $(\mathfrak g, [ \ ,
\ ])$ with a representation $(V, \rho^{L},\rho^{R})$ if $N$ is a
Nijenhuis operator and for every $y \in \mathfrak g$ and $w \in V$,
$(N, S)$ satisfies:
\begin{align} \label{3.3}
\rho^{L}(N(y)) S(w) &=S(\rho^{L}(N(y))
 w)+\rho^{L}(y)S^{2}(w)-S(\rho^{L}(y)S(w)),\\ \label{3.4}
\rho^{R}(N(y)) S(w) &=S(\rho^{R}(N(y))
 w)+\rho^{R}(y)S^{2}(w)-S(\rho^{R}(y)S(w)).
 \end{align}

Obviously, a trivial infinitesimal deformation gives rise to a
Nijenhuis pair. Conversely, a Nijenhuis pair also determines a
trivial infinitesimal deformation via
\eqref{deform1}-\eqref{deform3}. More precisely, we have the
following result.

{\bf Theorem 3.4.} Let $(N, S)$ be a Nijenhuis pair on a Leibniz
algebra $(\mathfrak g, [ \ , \ ])$ with a representation
 $(V,
\rho^{L},\rho^{R})$. Then the following defines a trivial
infinitesimal deformation of $(\mathfrak g, [ \ , \ ])$ with $(V,
\rho^{L},\rho^{R})$:
$$\omega(y_1,y_2)=[N(y_1), y_2]+[y_1, N(y_2)]-N[y_1, y_2], ~~\forall ~y_1,y_2\in \mathfrak g,$$ 
$$\varpi^{L}(y)=\rho^{L}(N(y))+\rho^{L}(y)S-S\rho^{L}(y)
, ~~\forall ~y \in\mathfrak g,$$
$$\varpi^{R}(y)=\rho^{R}(N(y))+\rho^{R}(y)S-S\rho^{R}(y),~~\forall ~y \in\mathfrak g
.$$

{\bf Theorem 3.5.} $(N, S)$ is a Nijenhuis pair on a Leibniz algebra
$(\mathfrak g, [ \ , \ ])$ with a representation
$(V,\rho^{L},\rho^{R})$ if and only if $(N, S^{*})$ is a
dual-Nijenhuis pair on $(\mathfrak g, [ \ , \ ])$ with
 the dual representation
$(V^{*},(\rho^{L})^{*},-(\rho^{L})^{*}-(\rho^{R})^{*})$.
\begin{proof}
If $(N, S)$ is a Nijenhuis pair, for every $y\in \mathfrak g,w\in
V,\alpha\in V^{*}$,
\begin{eqnarray*}&&\langle \rho^{L}(N(y)) S(w)
 -S(\rho^{L}(N(y))
 w)-S(\rho^{L}(y)S(w))+S^{2}(\rho^{L}(y)(w)),\alpha\rangle
 \\&=&\langle w,-S^{*}(\rho^{L})^{*}(N(y))(\alpha)
 +(\rho^{L})^{*}(N(y))(S^{*}(\alpha))+S^{*}(\rho^{L})^{*}(y)(S^{*}(\alpha))
 -(\rho^{L})^{*}(y)((S^{*})^{2}(\alpha))\rangle
 \\&=&0,
\end{eqnarray*}
 that is,
 $$S^{*}(\rho^{L})^{*}(N(y))(\alpha)
-(\rho^{L})^{*}(N(y))(S^{*}(\alpha))-S^{*}(\rho^{L})^{*}(y)(S^{*}(\alpha))+(\rho^{L})^{*}(y)((S^{*})^{2}(\alpha))=0.$$
Similarly, $$S^{*}(\rho^{R})^{*}(N(y))(\alpha)
-(\rho^{R})^{*}(N(y))(S^{*}(\alpha))-S^{*}(\rho^{R})^{*}(y)(S^{*}(\alpha))+(\rho^{R})^{*}(y)((S^{*})^{2}(\alpha))=0.$$
It follows that
\begin{eqnarray*}&&S^{*}((\rho^{R})^{*}+(\rho^{L})^{*})(N(y))(\alpha)
-((\rho^{R})^{*}+(\rho^{L})^{*})(N(y))(S^{*}(\alpha))-S^{*}((\rho^{R})^{*}+(\rho^{L})^{*})(y)(S^{*}(\alpha))
\\&&+((\rho^{R})^{*}+(\rho^{L})^{*})(y)((S^{*})^{2}(\alpha))=0. \end{eqnarray*}
Hence, $(N, S^{*})$ is a dual-Nijenhuis pair on $(\mathfrak g, [ \ ,
\ ])$ with
 the dual representation
$(V^{*},(\rho^{L})^{*},-(\rho^{L})^{*}-(\rho^{R})^{*})$. The converse
direction can be obtained similarly.
\end{proof}


{\bf Definition 3.6.} A {\it perfect Nijenhuis pair} on a Leibniz
algebra $(\mathfrak g, [ \ , \ ])$ with a representation
$(V,\rho^{L},\rho^{R})$ is a Nijenhuis pair $(N, S)$ with the
following identities:
$$S^{2}(\rho^{L}(y)(w))+\rho^{L}(y)S^{2}(w)=2S(\rho^{L}(y)S(w)),$$
$$S^{2}(\rho^{R}(y)w)+\rho^{R}(y)S^{2}(w)=2S(\rho^{R}(y)S(w)),~~\forall ~y \in\mathfrak g,w\in V.$$

 {\bf Proposition 3.7}. (i) Let $((\mathfrak g_{1}, [ \ , \ ]_1),(\mathfrak
g_{2},[ \ , \
]_2),\rho^{L}_{1},\rho^{R}_{1},\rho^{L}_{2},\rho^{R}_{2})$ be
 a matched pair of Leibniz algebras. Suppose that
$(N,S)$ is a Nijenhuis pair on $(\mathfrak g_{1}, [ \ , \ ]_1)$ with
the representation $(\mathfrak g_{2},\rho^{L}_{1},\rho^{R}_{1})$,
and $(S,N)$ a Nijenhuis pair on $(\mathfrak g_{2},[ \ , \ ]_2)$ with
the representation $(\mathfrak g_{1},\rho^{L}_{2},\rho^{R}_{2})$.
Then $N+S$ is a Nijenhuis operator on $\mathfrak g_{1}\bowtie
\mathfrak g_{2}$;

(ii) If $(N, S)$ is a perfect Nijenhuis pair on a
Leibniz algebra $(\mathfrak g, [ \ , \ ])$ with a representation
$(V,\rho^{L},\rho^{R})$, then $N + S^{*}$ is
 a Nijenhuis operator on the semidirect product
$\mathfrak g\ltimes V^{*}$, where $V^{*}$ is the dual representation
of $V$.
\begin{proof} 
We only check case (i), as case (ii) can be checked similarly. For
any $x_0,x_1\in \mathfrak g_{1},~a_0,a_1\in \mathfrak g_{2}$,
 \begin{eqnarray}\nonumber
&&(N+S)[x_0+a_0,x_1+a_1]_{N+S}\\ \nonumber
&=&(N+S)([N(x_0)+S(a_0),x_1+a_1]+[x_0+a_0, N(x_1)+S(a_1)]-(N+S)[x_0+a_0,x_1+a_1])\\
\nonumber &=&N[N(x_0),
x_1]_1+N(\rho^{L}_{2}(S(a_0))x_1+\rho^{R}_{2}(a_1)N(x_0))+S[S(a_0),a_1]_2+S(\rho^{L}_{1}(N(x_0))a_1
\\ \nonumber &&+\rho^{R}_{1}(x_1)S(a_0))+N[x_0,
N(x_1)]_1+N(\rho^{L}_{2}(a_0)N(x_1)+\rho^{R}_{2}(S(a_1))x_0)+S[a_0,S(a_1)]_2\\
\nonumber &&+S(\rho^{L}_{1}(x_0)S(a_1) +\rho^{R}_{1}(N(x_1))a_0) -
N^{2}[x_0,x_1]_1-N^{2}(\rho^{L}_{2}(a_0)x_1+\rho^{R}_{2}(a_1)x_0)\\&&-S^{2}[a_0,a_1]_2-S^{2}(\rho^{L}_{1}(x_0)a_1
+\rho^{R}_{1}(x_1)a_0).   \label{3.5} 
\end{eqnarray}
At the same time,
\begin{eqnarray} \nonumber
&&[(N+S)(x_0+a_0), (N+S)(x_1+a_1)] \\ \nonumber &=&[N(x_0),
N(x_1)]_1+\rho^{L}_{2}(S(a_0))N(x_1)+\rho^{R}_{2}(S(a_1))N(x_0)+[S(a_0),S(a_1)]_2\\
\label{3.6} &&+\rho^{L}_{1}(N(x_0))S(a_1)
+\rho^{R}_{1}(N(x_1))S(a_0).      
\end{eqnarray}
Combining \eqref{3.5} and \eqref{3.6}, we get
$$(N+S)[(x_0+a_0),(x_1+a_1)]_{N+S}-[(N+S)(x_0+a_0),(N+S)(x_1+a_1)]=0.$$
Hence $N+S$ is a Nijenhuis operator on $\mathfrak g_{1}\bowtie
\mathfrak g_{2}$.
\end{proof}

Let's consider the linear maps
$\hat{\rho^{L}},~\hat{\rho^{R}},~\tilde{\rho^{L}},~\tilde{\rho^{R}}:\mathfrak
g\longrightarrow \mathfrak{gl}(V)$ defined by: for every $y\in \mathfrak g$
\begin{equation}\label{e:4maps}
\begin{aligned}\hat{\rho^{L}}(y)&=\rho^{L}(N(y))+\rho^{L}(y)S-S\rho^{L}(y),\\
\hat{\rho^{R}}(y)&=\rho^{R}(N(y))+\rho^{R}(y)S-S\rho^{R}(y),\\
\tilde{\rho^{L}}(y)&=\rho^{L}(N(y))+S\rho^{L}(y)-\rho^{L}(y)S,\\
\tilde{\rho^{R}}(y)&=\rho^{R}(N(y))+S\rho^{R}(y)-\rho^{R}(y)S.
\end{aligned}
\end{equation}

{\bf Corollary 3.8.} (i) Let $(N, S)$ be a Nijenhuis
 pair on a Leibniz algebra $(\mathfrak g,[ \ , \ ])$ with a representation
 $(V,\rho^{L},\rho^{R})$, then $(V,\hat{\rho^{L}},\hat{\rho^{R}})$
   becomes a representation of the
Leibniz
 algebra $(\mathfrak g,[ \ , \ ]_{N})$ (cf.\eqref{e:Nijen});

(ii) Let $(N, S)$ be a dual-Nijenhuis
 pair on a Leibniz algebra $(\mathfrak g,[ \ , \ ])$ with a representation
 $(V,\rho^{L},\rho^{R})$, then
 $(V,\tilde{\rho^{L}},\tilde{\rho^{R}})$ becomes a representation of the
Leibniz
 algebra $(\mathfrak g,[ \ , \ ]_{N})$.
\begin{proof} (i) On the Nijenhuis pair $(N,S)$, for any $x_0,x_1\in \mathfrak
g,v_0,v_1\in V$,
\begin{eqnarray*}&&[x_0+v_0,x_1+v_1]_{N+S}\\&=&[(N+S)(x_0+v_0),x_1+v_1]+[x_0+v_0,(N+S)(x_1+v_1)]-(N+S)[x_0+v_0,x_1+v_1]
\\&=&[N(x_0),x_1]+\rho^{L}(N(x_0))v_1+\rho^{R}(x_1)S(v_0)+[x_0,N(x_1)]+\rho^{L}(x_0)S(v_1)\\&&+\rho^{R}(N(x_1))v_0
-N[x_0,x_1]-S(\rho^{L}(x_0)v_1+\rho^{R}(x_1)v_0)
\\&=&
[x_0,x_1]_{N}+\hat{\rho^{L}}(x_0)v_1+\hat{\rho^{R}}(x_1)v_0,
\end{eqnarray*} which indicates that
$(V,\hat{\rho^{L}},\hat{\rho^{R}})$ is a representation of
$(\mathfrak g,[ \ , \ ]_{N})$.

(ii) Clearly, the dual maps
$\tilde{\rho^{L}}^{*},\tilde{\rho^{R}}^{*}$
 of $\tilde{\rho^{L}},\tilde{\rho^{R}}$ are respectively given by
$$\tilde{\rho^{L}}^{*}(y)=(\rho^{L})^{*}(N(y))+(\rho^{L})^{*}(y)S^{*}-S^{*}(\rho^{L})^{*}(y),~~\forall~y\in\mathfrak g,$$ and $$
\tilde{\rho^{R}}^{*}(y)=(\rho^{R})^{*}(N(y))+(\rho^{R})^{*}(y)S^{*}-S^{*}(\rho^{R})^{*}(y).$$
Since $(N,S)$ is a dual-Nijenhuis pair on $(\mathfrak g,[ \ , \ ])$
with the representation $(V,\rho^{L},\rho^{R})$, according to
Theorem 3.5, $(N, S^{*})$ is a Nijenhuis pair on $(\mathfrak g,[ \ ,
\ ])$ with the representation
$(V^{*},(\rho^{L})^{*},-(\rho^{L})^{*}-(\rho^{R})^{*})$. In the
light of (i), we get the conclusion.
\end{proof}
The above statements show that though these results are natural
generalizations of the corresponding result in Lie algebras and
associative algebras, it is still worthy to carefully state them as
there are clear
differences in formulation of dual representations for Leibniz algebras.

\section{(Dual) $\mathcal{O}$N-structures and
Compatible $\mathcal{O}$-operators}

Let $K$ be an $\mathcal{O}$-operator on a Leibniz algebra
$(\mathfrak g, [ \ , \ ])$ associated to a representation
$(V,\rho^{L},\rho^{R})$. On the Leibniz algebra $(V_K, [ \ , \
]^K)$, let $S\in \mathfrak{gl}(V)$ and define
a deformed operation called $[ \ , \ ]_{S}^{K}$ on $V_{K}$ as follows: 
$$[w,u]_{S}^{K}=[S(w),u]^{K}+[w,S(u)]^{K}-S[w,u]^{K},~~\forall ~ w,u\in V.$$

Now let $(N, S)$ be a Nijenhuis pair or a dual Nijenhuis pair on the
Leibniz algebra $(\mathfrak g, [ \ , \ ])$ with the representation
$(V,\rho^{L},\rho^{R})$. Based on Corollary 3.8, there are two
representations $(V,\hat{\rho^{L}},\hat{\rho^{R}})$ and
$(V,\tilde{\rho^{L}},\tilde{\rho^{R}})$ of the deformed Leibniz
algebra $(\mathfrak g,[\ , \ ]_{N})$. We can define two operations
$[ \ , \ ]_{\hat{\rho}}^{K},[\ , \ ]_{\tilde{\rho}}^{K}:V\times
V\longrightarrow V$ by
\begin{equation}\label{e:rohhat}
\begin{aligned}
&[w,u]_{\hat{\rho}}^{K}=\hat{\rho^{L}}(K(w))(u)+\hat{\rho^{R}}(K(u))(w),\\
&[w,u]_{\tilde{\rho}}^{K}=\tilde{\rho^{L}}(K(w))(u)+\tilde{\rho^{R}}(K(u))(w),~~~\forall
~ w,u\in V.
\end{aligned}
\end{equation}

In general, $[ \ , \ ]_{S}^{K},[ \ , \ ]_{\hat{\rho}}^{K}$ and $[ \
, \ ]_{\tilde{\rho}}^{K}$ are not Leibniz brackets on $V$.
Obviously, if $S$ is
 a Nijenhuis operator on $(V_K, [ \ , \
]^K)$, then  $[ \ , \ ]_{S}^{K}$ is a Leibniz
 bracket. At the same time, $[ \ , \ ]_{\hat{\rho}}^{K}$ and $[ \ , \
]_{\tilde{\rho}}^{K}$ become Leibniz
 brackets when $K$ is an $\mathcal{O}$-operator on $(\mathfrak g, [ \ , \ ]_{N})$ associated to
$(V,\hat{\rho^{L}},\hat{\rho^{R}})$ and
$(V,\tilde{\rho^{L}},\tilde{\rho^{R}})$ respectively. In the
following, we will study when this happens. 

{\bf Definition 4.1.} Let $K:V\longrightarrow \mathfrak g$ be an
$\mathcal{O}$-operator on a Leibniz algebra $(\mathfrak g, [ \ , \
])$ associated to a representation $(V,\rho^{L},\rho^{R})$ and $(N,
S)$ a Nijenhuis pair (resp. a dual Nijenhuis pair). The triple $(K,
S, N)$ is called an {\it $\mathcal{O}$N-structure (resp. a dual
$\mathcal{O}$N-structure)} on $(\mathfrak g, [ \ , \ ])$ associated
to $(V,\rho^{L},\rho^{R})$ if the following holds:
\begin{equation}\label{4.1}
NK=KS ~~\hbox{and }~~[w,u]^{NK}=[w,u]_{S}^{K},   ~~~\forall ~ w,u\in
V,
\end{equation}
where
$[w,u]^{NK}=\rho^{L}(NK(w))(u)+\rho^{R}(NK(u))(w).$

In the definition of an $\mathcal{O}$N-structure (resp. a dual
$\mathcal{O}$N-structure) $(K,S,N)$ on $(\mathfrak g, [ \ , \ ])$
associated to a representation $(V,\rho^{L},\rho^{R})$, $NK$ does
not need to be assumed a priori an $\mathcal{O}$-operator on
$(\mathfrak g, [ \ , \ ])$ associated to the representation
$(V,\rho^{L},\rho^{R})$. This is due to the fact that when $(K,S,N)$
is an $\mathcal{O}$N-structure or a dual $\mathcal{O}$N-structure,
$NK$ will be an $\mathcal{O}$-operator on $(\mathfrak g, [ \ , \ ])$
associated to $(V,\rho^{L},\rho^{R})$. We will discuss this in
Proposition 4.4.

 Similar to Lie algebras [HLS], we can show that

{\bf Proposition 4.2.} (i) If $(K, S, N)$ is an
$\mathcal{O}$N-structure on a Leibniz algebra $(\mathfrak g, [ \ , \
])$ associated to a representation $(V,\rho^{L},\rho^{R})$,
 then
$$[w,u]_{S}^{K}=[w,u]_{\hat{\rho}}^{K},~~~~\forall ~
w,u\in V;$$

(ii) If $(K, S, N)$ is a dual $\mathcal{O}$N-structure on a Leibniz
algebra $(\mathfrak g, [ \ , \ ])$ associated to a representation
$(V,\rho^{L},\rho^{R})$, then
$$[w,u]_{S}^{K}=[w,u]_{\tilde{\rho}}^{K},~~~\forall ~
w,u\in V.$$

Thus the operations $[ \ , \
]_{\hat{\rho}}^{K}$ = $[ \ , \ ]_{S}^{K}$ when $(K, S, N)$ is an
$\mathcal{O}$N-structure and $[ \ , \ ]_{S}^{K}=[ \ , \
]_{\tilde{\rho}}^{K}$ if $(K, S, N)$ is a dual
$\mathcal{O}$N-structure (cf.\eqref{e:rohhat}). Moreover, these operations are Leibniz
brackets.

{\bf Theorem 4.3} If $(K, S, N)$ is an $\mathcal{O}$N-structure or a
dual $\mathcal{O}$N-structure on a Leibniz algebra $(\mathfrak g, [
\ , \ ])$ associated to a representation $(V,\rho^{L},\rho^{R})$,
then $S$ is a Nijenhuis operator on $(V_{K},[ \ , \ ]^{K})$.
Moreover, the operations $[ \ , \ ]_{S}^{K}, [ \ , \
]_{\tilde{\rho}}^{K}$ and $[ \ , \ ]_{\hat{\rho}}^{K}$ all give rise
to Leibniz algebras.
\begin{proof} When $(K, S, N)$ is an $\mathcal{O}$N-structure,
substituting $y$ by $K(v)$ in \eqref{3.1}, we have that
\begin{eqnarray*}&&\rho^{L}(NK(v))S(w)-S(\rho^{L}(NK(v))w)-S(\rho^{L}(K(v))S(w))+S^{2}(\rho^{L}(K(v))w)\\&=&
\rho^{L}(KS(v))S(w)-S(\rho^{L}(KS(v))w)-S(\rho^{L}(K(v))S(w))+S^{2}(\rho^{L}(K(v))w)\\
&=&S(v)\lhd^{K}S(w)-S(S(v)\lhd ^{K}w)-S(v\lhd^{K}
S(w))+S^{2}(v\lhd^{K} w). 
\end{eqnarray*}
The same identity holds by replacing $\lhd^K$ with $\rhd^K$. Since
$[ \ , \ ]^{K}=\lhd^{K}+\rhd^{K}$ and \eqref{3.1},
$$[S(v),S(w)]^{K}=S([S(v),w]^{K} )+S([v,S(w)]^{K})
-S^{2}([v,w]^{K} ).$$ Thus $S$ is a Nijenhuis operator on $(V_{K}, [
\ ,\ ]^{ K })$.

If $(K, S, N)$ is a dual $\mathcal{O}$N-structure,
 it follows from $[v,w]^{NK}=[v,w]_{S}^{K} $ that
\begin{equation} \label{4.2}
\rho^{R}(K(w))S(v)+\rho^{L}(K(v))S(w)-S(\rho^{R}(K(w))v)-S(\rho^{L}(K(v))w)=0.
\end{equation}
Replacing $v$ by $S(v)$ in \eqref{4.2}, we have
\begin{equation}\label{4.3}
\rho^{R}(K(w))S^{2}(v)+\rho^{L}(K(S(v)))S(w)-S(\rho^{R}(K(w))S(v))-S(\rho^{L}(KS(v))w)=0.
\end{equation}
Applying $S$ to both sides of \eqref{4.2}, we get
\begin{equation}\label{4.4}
S(\rho^{R}(K(w))S(v))+S(\rho^{L}(K(v))S(w))-S^{2}(\rho^{R}(K(w))v)-S^{2}(\rho^{L}(K(v))w)=0.
\end{equation}
Replacing $y$ by $K(w)$
and $w$ by $v$ in the identities \eqref{3.3} and \eqref{3.4}, we get that
\begin{align}\label{4.5}
\rho^{L}(KS(w)) S(v)
 &=S(\rho^{L}(KS(w))
 v)+\rho^{L}(K(w))S^{2}(v)-S(\rho^{L}(K(w))S(v)),\\ \label{4.6}
\rho^{R}(KS(w)) S(v)
 &=S(\rho^{R}(KS(w))
 v)+\rho^{R}(K(w))S^{2}(v)-S(\rho^{R}(K(w))S(v)).
\end{align}
In view of \eqref{4.3}-\eqref{4.6},
\begin{eqnarray*}
&&[S(v),S(w)]^{K}-S([v,w]_{S}^{K})
\\&=&\rho^{L}(KS(v))S(w)+\rho^{R}(KS(w))S(v)-S(\rho^{L}(KS(v))w+\rho^{R}(K(w))S(v)\\&&+\rho^{L}(K(v))S(w)
+\rho^{R}(KS(w))v)+S^{2}(\rho^{L}(K(v))w+\rho^{R}(K(w))v)\\
&=&\rho^{L}(KS(v))S(w)+\rho^{R}(KS(w))S(v)-\rho^{R}(K(w))S^{2}(v)-\rho^{L}(KS(v))S(w)
\\&&-S(\rho^{L}(K(v))S(w) +\rho^{R}(KS(w))v)+
S(\rho^{R}(K(w))S(v))+S(\rho^{L}(K(v))S(w))\\&=&0.
\end{eqnarray*}
Thus the result follows.
\end{proof}

{\bf Proposition  4.4}. Suppose that $(K, S, N)$ is an
$\mathcal{O}$N-structure (resp. a dual $\mathcal{O}$N-structure) on
a Leibniz algebra $(\mathfrak g, [ \ , \ ])$ associated with a
representation $(V,\rho^{L},\rho^{R})$. Then

(i) $K$ is an $\mathcal{O}$-operator on $(\mathfrak g, [ \ , \
]_{N})$ associated to the representation
$(V,\hat{\rho^{L}},\hat{\rho^{R}})$ (resp. $(V,\tilde{\rho^{L}},\tilde{\rho^{R}})$) (cf. \eqref{e:4maps});


(ii) $NK$ is an $\mathcal{O}$-operator  $(\mathfrak g, [ \ , \ ])$
with respect to the representation $(V,\rho^{L},\rho^{R})$.

\begin{proof} These hold by Proposition 4.2 and direct calculation. \end{proof}

 {\bf Theorem 4.5.} Suppose that $(K, S, N)$ is an
$\mathcal{O}$N-structure on a Leibniz algebra $(\mathfrak g, [ \ , \
])$ with respect to a representation $(V,\rho^{L},\rho^{R})$ with
invertible $K$, then $(K, S, N)$ becomes a dual
$\mathcal{O}$N-structure.

\begin{proof}
It is enough to show that $(N, S)$ is a dual-Nijenhuis pair. Thanks
to the $\mathcal{O}$N-structure $(K, S, N)$, substituting $S(v)$ for
$v$ in \eqref{4.2} we have that
\begin{equation}\label{4.7}\rho^{R}(K(w))S^{2}(v)+\rho^{L}(KS(v))S(w)=S(\rho^{L}(KS(v))w+\rho^{R}(K(w))S(v)).\end{equation}
Taking account of $[v,w]_{S}^{K}=[v,w]^{KS}$,
$$S([v,w]^{KS})=S([v,w]_{S}^{K})=[S(v),S(w)]^{K},$$
that is,
\begin{equation}\label{4.8}
S(\rho^{L}(KS(v))w+\rho^{R}(KS(w))v)=\rho^{L}(KS(v))S(w)+\rho^{R}(KS(w))S(v).
\end{equation}
Combining \eqref{4.7} and \eqref{4.8}, we have
\begin{eqnarray}\nonumber
0&=&\rho^{R}(K(w))S^{2}(v)+S(\rho^{R}(KS(w))v)-\rho^{R}(KS(w))S(v)-S(\rho^{R}(K(w))S(v))\\
\label{4.9} &=&\rho^{R}(K(w))S^{2}(v)+S(\rho^{R}(NK(w))v)
-\rho^{R}(NK(w))S(v)-S(\rho^{R}(K(w))S(v)).
\end{eqnarray}
Since $K$ is invertible, letting $K(w)=x$, we get that
$$\rho^{R}(x)S^{2}(v)-\rho^{R}(N(x))S(v)+S(\rho^{R}(N(x))v)-S(\rho^{R}(x)S(v))=0.$$
 Similarly, the same identity holds for $\rho^L$,
thus $(N,S)$ is a dual-Nijenhuis pair.
\end{proof}

 {\bf Definition 4.6.} Let $K_1, K_2 $ be two $\mathcal{O}$-operators on a
Leibniz algebra $(\mathfrak g, [ \ , \ ])$ with respect to a
representation $(V,\rho^{L},\rho^{R})$. We say that $K_1$ and $K_2$
are {\it compatible} if
 $n_1K_1+n_2K_2$ is an $\mathcal{O}$-operator for any $n_1,
 n_2\in k$.

Similar to Lie algebras and associative algebras, the following results
are clear.

{\bf Proposition 4.7.} Suppose $K$ is an $\mathcal{O}$-operator on a
Leibniz algebra $(\mathfrak g, [ \ , \ ])$ with respect to a
representation $(V,\rho^{L},\rho^{R})$ and $S\in \mathfrak{gl}(V)$.
Then $K$ and $KS$ are compatible $\mathcal{O}$-operators if $(K, S,
N)$ is an $\mathcal{O}$N-structure or a dual
$\mathcal{O}$N-structure.

 {\bf Proposition 4.8.} Let $K_1, K_2$
 be two $\mathcal{O}$-operators on a Leibniz algebra $(\mathfrak g, [ \ , \
])$ with respect to a representation $(V,\rho^{L},\rho^{R})$. If
$K_1$ and $K_2$ are compatible and $K_1$ is invertible, then

(i) $(K_1,S=K_1^{-1}K_2,N=K_2K_1^{-1}) $ is a dual
$\mathcal{O}$N-structure;

(ii) $(K_2,S=K_1^{-1}K_2,N=K_2K_1^{-1})$ is a dual
$\mathcal{O}$N-structure.

\section{Strong Maurer-Cartan equation and dual $\mathcal{O}$N-structures}

Unless otherwise specified, we follow the same notations in Section 2.

 Let $(\mathfrak g_{1}, [\ , \ ]_1, \mathfrak g_{2}, [\ , \ ]_2, \rho^{L}_{1},\rho^{R}_{1},\rho^{L}_{2},\rho^{R}_{2})$ be
a matched pair of Leibniz algebras (or equivalently, twilled Leibniz
algebra $(\mathfrak g_{1}\bowtie\mathfrak g_{2},\mathfrak g_{1},
\mathfrak g_{2})$). It follows that $ C^{*}(\mathfrak g_1\bowtie
\mathfrak g_2,\mathfrak g_1\bowtie \mathfrak g_2)$ is a
 dgLa. Let $\vartheta=\hat{\mu}_1+\hat{\mu}_2\in C^{2}(\mathfrak g_1\bowtie \mathfrak g_2,\mathfrak g_1\bowtie \mathfrak g_2)$,
where $\hat{\mu}_1,\hat{\mu}_2$
are the Leibniz brackets on semidirect products $\mathfrak
g_1\ltimes \mathfrak g_2$
 and $\mathfrak g_2\ltimes \mathfrak g_1$ respectively.

Consider the graded vector space $C^{*}(\mathfrak g_1,\mathfrak
g_2)=\oplus_{p=1}^{\infty} C^{p}(\mathfrak g_1,\mathfrak
 g_2)=\oplus_{p=1}^{\infty}\hbox{Hom}(\otimes^{p}\mathfrak g_1,\mathfrak
 g_2).$ In view of Theorem 2.2, $(C^{*}(\mathfrak g_1,\mathfrak g_2), \{\{
\ , \}\}, d_{\hat{\mu}_1})$ is a dgLa with
 $$\{\{\varphi,\psi\}\}=(-1)^{m}\{\{\hat{\mu}_2,\hat{\varphi}\}^{B},\hat{\psi}\}^{B}$$
 for $\varphi\in C^{m}(\mathfrak g_1,\mathfrak
 g_2),\psi\in C^{n}(\mathfrak g_1,\mathfrak
 g_2)$. Here the differential $d_{\hat{\mu}_1}: C^{m}(\mathfrak g_1,\mathfrak
 g_2)\longrightarrow  C^{m+1}(\mathfrak g_1,\mathfrak
 g_2)$ is given by
 $d_{\hat{\mu}_1}(\varphi)=\{\hat{\mu}_1,\hat{\varphi}\}^{B}$.

Let $\Theta:\mathfrak g_{1}\longrightarrow \mathfrak g_{2}$ be a
linear map. The
 equations
 $$d_{\hat{\mu}_1}\Theta+\frac{1}{2}\{\{\Theta,\Theta\}\}=0,~~d_{\hat{\mu}_1}\Theta=\frac{1}{2}\{\{\Theta,\Theta\}\}=0$$
 are called {\it the Maurer-Cartan equation} and {\it the strong Maurer-Cartan
 equation } respectively.

{\bf Proposition 5.1.} Under the same notations above, we have

 (i) $\Theta$ is a solution of the Maurer-Cartan equation if and only if
 $\Theta$ satisfies
\begin{equation}\label{2.1}
[\Theta(y),\Theta(z)]_{2}+\rho_{1}^{L}(y)\Theta(z)+\rho_{1}^{R}(z)\Theta(y)
=\Theta(\rho_{2}^{L}\Theta(y)(z)+\rho_{2}^{R}\Theta(z)(y))+\Theta([y,z]_{1})
\end{equation}
for any $y,z\in \mathfrak g_{1}$;

(ii) $\Theta$ is a solution of the strong Maurer-Cartan equation if
and only if
 Eq. \eqref{2.1} holds and for any $y,z\in \mathfrak g_{1}$, the following condition holds:
\begin{equation}\label{2.2}
\Theta([y,z]_{1})=\rho_{1}^{L}(y)\Theta(z)+\rho_{1}^{R}(z)\Theta(y).
\end{equation}

\begin{proof} For any $y,z\in \mathfrak g_{1},a,b\in \mathfrak g_{2},$
\begin{eqnarray}\nonumber &&d_{\hat{\mu}_1}\Theta(y+a,z+b)=\{\hat{\mu}_1,\hat{\Theta}\}^{B}(y+a,z+b)\nonumber\\&=&(\hat{\mu}_1(\hat{\Theta}\otimes
I)+\hat{\mu}_1(I\otimes
\hat{\Theta})-\hat{\Theta}\hat{\mu}_1)(y+a,z+b)\nonumber\\&=&\rho_{1}^{R}(z)\Theta(y)+\rho_{1}^{L}(y)\Theta(z)-\Theta([y,z]_1)\nonumber,
\end{eqnarray}
that is, $d_{\hat{\mu}_1}\Theta=\frac{1}{2}\{\{\Theta,\Theta\}\}=0$
is equivalent to that (5.2) holds. Analogously,
$d_{\hat{\mu}_1}\Theta+\frac{1}{2}\{\{\Theta,\Theta\}\}=0$ is
equivalent to that (5.1) holds.
\end{proof}

 The following result is a special case of Proposition 5.1
 with $\mathfrak g_1=\mathfrak g$ and $\mathfrak
g_2=V_K$.

{\bf Corollary 5.2.} Let $K$ be an $\mathcal{O}$-operator on a
Leibniz algebra $\mathfrak g$ with a representation $V$, and let $
\Theta: \mathfrak g \longrightarrow V$ be a linear map. Then
$\Theta$ satisfies the strong Maurer-Cartan equation on the twilled
algebra $\mathfrak g\bowtie V_{K}$ if and only if
\begin{align}\label{2.3}
\Theta([y,z])&=\rho^{L}(y)\Theta(z)+\rho^{R}(z)\Theta(y),\\
\label{2.4}
[\Theta(y),\Theta(z)]^{K}&=\Theta(\varrho_{K}^{L}(\Theta(y))z+\varrho_{K}^{R}(\Theta(z))y),~~\forall~
y,z\in \mathfrak g.
\end{align}

According to \eqref{2.4}, a solution $\Theta: \mathfrak g \longrightarrow
V$ of the strong Maurer-Cartan equation on the twilled Leibniz
algebra $\mathfrak g\bowtie V_{K}$ becomes an
 $\mathcal{O}$-operator on the Leibniz
algebra $(V_{K}, [ \ , \ ]^K )$ with respect to the representation
$(\mathfrak g, \varrho_{K}^{L},\varrho_{K}^{R})$. The
$\mathcal{O}$-operator $\Theta$ leads to a Leibniz algebra structure
on $\mathfrak g$ with
$[y,z]^{\Theta}=\varrho_{K}^{L}(\Theta(y))z+\varrho_{K}^{R}(\Theta(z))y~(\forall
~y,z\in \mathfrak g)$. Denote this Leibniz algebra by $(\mathfrak
g_{\Theta},[ \ , \ ]^{\Theta})$.

Similar to $V_{K}$, $\rho_{\Theta}^{L},\rho_{\Theta}^{R}:\mathfrak
g_{\Theta}\longrightarrow \mathfrak{gl}(V_{K})$ afford a
representation of $\mathfrak g_{\Theta}$ on $V_{K}$ with
$$ \rho_{\Theta}^{L}(y)(w)=[\Theta(y),w]^{K}-\Theta(\varrho_{K}^{R}(w)y),~
\rho_{\Theta}^{R}(y)(w)=[w,\Theta(y)]^{K}-\Theta(\varrho_{K}^{L}(w)y)$$
for any $y\in \mathfrak g,w\in V$. Then $(\mathfrak g\oplus V,\{ \ ,
\ \}_{K}^{\Theta})$ is a Leibniz algebra with
\begin{equation*}
\{w_0+x_0,w_1+x_1
\}_{K}^{\Theta}=[w_0,w_1]^{K}+\varrho_{K}^{L}(w_0)x_1+\varrho_{K}^{R}
(w_1)x_0+[x_0,x_1]^{\Theta}+\rho_{\Theta}^{L}(x_0)w_1+\rho_{\Theta}^{R}(x_1)w_0,
\end{equation*}
for every $x_0,x_1\in \mathfrak g,w_0,w_1\in V.$

 {\bf Proposition 5.3.} Suppose that $\Theta:\mathfrak
g\longrightarrow V$ is a solution of the strong Maurer-Cartan
equation on the twilled Leibniz algebra $(\mathfrak g\bowtie
V_{K},\{ \ , \ \}_{K}$). Then

(i) $(\mathfrak g\oplus V,\{ \ , \ \}_{K}^{\Theta})$ is a Leibniz
algebra, denoted by $( V_{K}\bowtie\mathfrak g_{\Theta},\{ \ ,
\ \}_{K}^{\Theta})$.

(ii) $K$ is a solution of the strong Maurer-Cartan equation on
$(V_{K}\bowtie\mathfrak g_{\Theta},\{ \ , \ \}_{K}^{\Theta})$.

(iii) $K$ is an $\mathcal{O}$-operator on the Leibniz algebra
$(\mathfrak g_{\Theta},[ \ , \ ]^{\Theta})$ associated to the
representation $(V_{K},\rho_{\Theta}^{L}, \rho_{\Theta}^{R})$.
\begin{proof} (i) was made clear earlier. Let's consider
(ii). Note that $K$ is an $\mathcal{O}$-operator on the Leibniz
algebra $(\mathfrak g,[ \ , \ ])$ associated to the representation
$(V,\rho^{L},\rho^{R})$, using \eqref{e:natural}, for any
$w_1,w_2\in V$,
\begin{eqnarray}
\nonumber&&\varrho_{K}^{L}(w_1)K(w_2)+\varrho_{K}^{R}(w_2)K(w_1)\nonumber
\\&=&[K(w_1),K(w_2)]-K(\rho^{R}(K(w_2))w_1)+[K(w_1),K(w_2)]-K(\rho^{L}(K(w_1))w_2)
\nonumber\label{6.1}\\ &=&[K(w_1),K(w_2)]=K[w_1,w_2]^{K}.
\end{eqnarray}
Since $\Theta:\mathfrak g\longrightarrow V$ is a solution of the
strong Maurer-Cartan equation on the twilled Leibniz algebra
$\mathfrak g\bowtie V_{K}$, it follows from \eqref{2.4} that
$\Theta$ is an $\mathcal{O}$-operator on the Leibniz algebra
$(V_{K}, [ \ , \ ]^K )$ with respect to the representation
$(\mathfrak g, \varrho_{K}^{L},\varrho_{K}^{R})$. Hence, in the
light of \eqref{e:natural} and \eqref{2.3}, for every $w_1,w_2\in
V$,
\begin{align}
&[K(w_1),K(w_2)]^{\Theta}\nonumber=\varrho_{K}^{L}(\Theta
K(w_1))K(w_2)+\varrho_{K}^{R}(\Theta K(w_2))K(w_1)
\nonumber\\
&=[K\Theta K(w_1),K(w_2)]-K(\rho^{R}(K(w_2))\Theta
K(w_1))+[K(w_1),K\Theta K(w_2)] \nonumber
\\
&\nonumber\qquad-K(\rho^{L}(K(w_1))\Theta K(w_2))\label{6.2}\\
&= [K\Theta
K(w_1),K(w_2)]+[K(w_1),K\Theta K(w_2)]-K\Theta[K(w_1),K(w_2)].
\end{align}
Meanwhile, $\Theta$ is an $\mathcal{O}$-operator on the Leibniz
algebra $(V_{K},[ \ , \ ]^{K})$ with respect to the representation
$(\mathfrak g,\varrho^{L}_{K},\varrho^{R}_{K})$, for every
$w_1,w_2\in V$,
\begin{eqnarray}\nonumber
&&K(\rho_{\Theta}^{L}(K(w_1))w_2+\rho_{\Theta}^{R}(K(w_2))w_1)\\
\nonumber &=&K\left([\Theta
K(w_1),w_2]^{K}-\Theta(\varrho_{K}^{R}(w_2)K(w_1))+[w_1,\Theta
K(w_2)]^{K}-\Theta(\varrho_{K}^{L}(w_1)K(w_2))\right) \\ \nonumber
&=&[K\Theta K(w_1),K(w_2)]+[K(w_1),K\Theta
K(w_2)]-K\Theta([K(w_1),K(w_2)]-K(\rho^{L}(K(w_1))w_2)\\ \nonumber
&&+[K(w_1),K(w_2)]-K(\rho^{R}(K(w_2))w_1)\label{6.3} \\
&=&[K\Theta K(w_1),K(w_2)]+[K(w_1),K\Theta K(w_2)]-K\Theta[K(w_1),K(w_2)].
\end{eqnarray}
Combining \eqref{6.2} and \eqref{6.3}, we have that
\begin{equation} \label{6.4}
[K(w_1),K(w_2)]^{\Theta}=K(\rho_{\Theta}^{L}(K(w_1))w_2+\rho_{\Theta}^{R}(K(w_2))w_1).
\end{equation}
Then \eqref{6.1} and \eqref{6.4} imply (ii) in view of Corollary
5.2. (iii) follows directly from \eqref{6.4}.
\end{proof}

{\bf Theorem 5.4.} Let $K$ be an $\mathcal{O}$-operator on a Leibniz
algebra $(\mathfrak g, [ \ , \ ])$ with respect to a representation
$(V, \rho^{L},\rho^{R})$, and $\Theta:\mathfrak g\longrightarrow V$
a solution of the strong Maurer-Cartan equation on the twilled
Leibniz algebra $(\mathfrak g\bowtie V_{K},\{ \ , \ \}_{K})$. Then

(i) $(K, S, N)$ is a dual $\mathcal{O}$N-structure on $(\mathfrak g,
[ \ , \ ])$ with respect to the representation $(V,
\rho^{L},\rho^{R})$, where $N = K\Theta$ and $S = \Theta K$.

(ii) $(\Theta, S, N)$ is a dual $\mathcal{O}$N-structure on $(V_K ,
[ \ , \ ]^{K} )$ associated to the representation $(\mathfrak g,
\varrho_{K}^{L} ,\varrho_{K}^{R} )$, where $N = K\Theta$ and $S=
\Theta K$.
\begin{proof} (i) Note that $\Theta$ is an $\mathcal{O}$-operator on $(V_{K}, [ \ , \ ]^K )$ with respect to the
representation $(\mathfrak g, \varrho_{K}^{L},\varrho_{K}^{R})$. It
follows from \eqref{2.3} and (5.4) that
\begin{eqnarray*}[K\Theta(y),K\Theta(z)]&=&K[\Theta(y),\Theta(z)]^{K}\\&=&K\Theta(\varrho_{K}^{L}(\Theta(y))z+\varrho_{K}^{R}(\Theta(z))y)
\\&=&K\Theta([K\Theta(y),z]-K(\rho^{R}(z)\Theta(y)))+K\Theta([y,K\Theta(z)]-K(\rho^{L}(y)\Theta(z)))
\\&=&K\Theta([K\Theta(y),z]+[y,K\Theta(z)]-K\Theta[y,z]),~~\forall
y,~z\in \mathfrak g,
\end{eqnarray*}
which implies that $K\Theta$ is a Nijenhuis operator on $(\mathfrak
g, [ \ , \ ])$.

We now check that  \eqref{3.3}, \eqref{3.4} and \eqref{4.1} hold for
$N = K\Theta$ and $S = \Theta K$. In fact, by \eqref{2.3} and $K$
being an $\mathcal{O}$-operator it follows that for every $w,u\in
V$,
\begin{eqnarray}\nonumber
\Theta K(\rho^{L}(K(w))u+\rho^{R}(K(u))w)&=& \Theta [K(w),K(u)]
\\ \label{6.5} 
&=&\rho^{L}(K(w))\Theta K(u)+\rho^{R}(K(u))\Theta K(w).
\end{eqnarray}
Replacing $y$ by $K(w)$ in (5.4) and using \eqref{2.3}, we have that
\begin{eqnarray*}&&[\Theta K(w),\Theta (z)]^{K}-\Theta(\varrho_{K}^{R}(\Theta (z))K(w)+\varrho_{K}^{L}(\Theta
K(w))z)
\\&=&
\rho^{R}(K\Theta (z))\Theta K(w)+\rho^{L}(K\Theta K(w))\Theta
(z)-\Theta ([K(w),K\Theta(z)]-K(\rho^{L}(K(w))\Theta(z)))\\&&-
\Theta([K\Theta K(w),z]-K(\rho^{R}(z)\Theta K(w)))
\\&=&\rho^{R}(K\Theta (z))\Theta K(w)+\rho^{L}(K\Theta K(w))\Theta
(z)-\rho^{L}(K(w))\Theta K\Theta (z)-\rho^{R}(K\Theta z)\Theta K(w)
\\&&+\Theta K(\rho^{L}(K(w))\Theta (z)) -\rho^{L}(K \Theta K(w))\Theta
(z)-\rho^{R}(z)\Theta K \Theta K(w) +\Theta K(\rho^{R}(z)\Theta
K(w))
\\&=&-\rho^{L}(K(w))\Theta K\Theta (z)+\Theta K(\rho^{L}(K(w))\Theta (z))-\rho^{R}(z)\Theta K \Theta K(w)
+\Theta K(\rho^{R}(z)\Theta K(w))=0.
\end{eqnarray*}
Therefore,
\begin{equation}\label{6.6a} 
-\rho^{L}(K(w))\Theta K\Theta (z)+\Theta K(\rho^{L}(K(w))\Theta
(z))=\rho^{R}(z)\Theta K \Theta K(w) -\Theta K(\rho^{R}(z)\Theta
K(w)).
\end{equation}
Combining \eqref{6.5} and \eqref{6.6a},
 we obtain that
\begin{equation}\label{6.7} 
\rho^{R}(K\Theta(z))\Theta K(w)-\Theta
K(\rho^{R}(K\Theta(z))w)=\rho^{R}(z)\Theta K\Theta K(w)-\Theta
K(\rho^{R}(z)\Theta K(w)).
\end{equation}
On the other hand, replacing $z$ with $K(u)$ in (5.4), it follows
from \eqref{2.3} that
\begin{equation}\label{6.8aa}          
-\rho^{R}(K(u))\Theta K\Theta(y)+\Theta
K(\rho^{R}(K(u))\Theta(y))=\rho^{L}(y)\Theta K\Theta K(u)-\Theta
K(\rho^{L}(y)\Theta K(u)). 
\end{equation}
Combining \eqref{6.5} and \eqref{6.8aa}, we find that
\begin{equation}\label{6.6}
\rho^{L}(K\Theta(y))\Theta K(u)-\Theta K(\rho^{L}(K\Theta(y))u)
=\rho^{L}(y)\Theta K\Theta K(u)-\Theta K(\rho^{L}(y)\Theta K(u)).
\end{equation}
Thus, \eqref{6.7} and \eqref{6.6} yield \eqref{3.3} and \eqref{3.4}
respectively.

At the same time,
\begin{eqnarray*}[w,u]^{K}_{S}-[w,u]^{KS}&=&\rho^{L}(KS(w))u+\rho^{R}(K(u))S(w)+\rho^{L}(K(w))S(u)+\rho^{R}(KS(u))w
\\&&-S\rho^{L}(K(w))u-S\rho^{R}(K(u))w-\rho^{L}(KS(w))u-\rho^{R}(KS(u))w
\\&=&\rho^{R}(K(u))S(w)+\rho^{L}(K(w))S(u)
-S\rho^{L}(K(w))u-S\rho^{R}(K(u))w\\&=& \rho^{L}(K(w))\Theta
K(u)+\rho^{R}(K(u))\Theta K(w)-\Theta
K(\rho^{L}(K(w))u+\rho^{R}(K(u))w)\\&=&0,\end{eqnarray*} where we
have used \eqref{6.5} in the last equation.

Similarly, (ii) can be verified.
\end{proof}

The following is expected from an analogous result for Lie algebras,
we include its proof as it needs explicit formulation.

 {\bf Theorem 5.5.} Suppose that $(K, S, N)$ is a dual
$\mathcal{O}$N-structure on a Leibniz algebra $(\mathfrak g, [ \ , \
])$ with respect to a representation $(V, \rho^{L},\rho^{R})$
 and $K$ is invertible, then
$\Theta=K^{-1}N=SK^{-1}:\mathfrak g\longrightarrow V$ is a solution
of the strong Maurer-Cartan equation on $\mathfrak g\bowtie V_{K}$.

\begin{proof} As $N = K \Theta$ is a Nijenhuis operator on $(\mathfrak g, [ \ , \
])$, for any $x_0,x_1\in \mathfrak g$,
\begin{equation}\label{6.8a}
[K\Theta(x_0),K\Theta(x_1)]=K\Theta([K\Theta(x_0),x_1]+[x_0,K\Theta(x_1)]-K\Theta[x_0,x_1]).
\end{equation}

Note that
$K$ is an invertible
$\mathcal{O}$-operator on $(\mathfrak g, [ \ , \ ])$ associated with
the representation $(V, \rho^{L},\rho^{R})$ and
$\Theta=K^{-1}N=SK^{-1}$. Let $v, w$ be the respective preimages  of $x_0, x_1$, then
\begin{eqnarray} \nonumber
K\Theta[x_0,x_1]&=&K\Theta[K(v),K(w)]=K\Theta
K(\rho^{L}(K(v))w+\rho^{R}(K(w))v)
\\ \nonumber &=&KS(\rho^{L}(K(v))w+\rho^{R}(K(w))v)
\\ \nonumber &=&K(\rho^{L}(K(v))S(w)+\rho^{R}(K(w))S(v))
\\\label{6.9a}
&=&K(\rho^{L}(x_0)\Theta(x_1)+\rho^{R}(x_1)\Theta(x_0)),
\end{eqnarray}
where we have used \eqref{4.2}. According
to \eqref{6.9a} and the definitions of $\varrho_{K}^{L}$ and
$\varrho_{K}^{R}$,
\begin{eqnarray} \nonumber
\varrho^{L}_{K}(\Theta(x_0))x_1+\varrho^{R}_{K}(\Theta(x_1))x_0&=&[K\Theta(x_0),x_1]-K(\rho^{R}(x_1)\Theta(x_0))
+[x_0,K\Theta(x_1)]\\ \nonumber &&-K(\rho^{L}(x_0)\Theta(x_1))
\\\label{6.9}
&=& [K\Theta(x_0),x_1]+[x_0,K\Theta(x_1)]-K\Theta[x_0,x_1].
\end{eqnarray}
Noting that $K$ is an $\mathcal{O}$-operator, \eqref{6.8a} and
\eqref{6.9} imply that
\begin{equation}\label{6.10a}K[\Theta(x_0),\Theta(x_1)]^{K}=[K\Theta(x_0),K\Theta(x_1)]
=K\Theta(\varrho^{L}_{K}(\Theta(x_0))x_1+\varrho^{R}_{K}(\Theta(x_1))x_0).\end{equation}
As $K$ is invertible, \eqref{6.9a} and \eqref{6.10a} indicate that
\begin{equation}\label{6.13}
\Theta[x_0,x_1]=\rho^{L}(x_0)\Theta(x_1)+\rho^{R}(x_1)\Theta(x_0),
\end{equation}
and
\begin{equation}\label{6.10}
[\Theta(x_0),\Theta(x_1)]^{K}=\Theta(\varrho^{L}_{K}(\Theta(x_0))x_1+\varrho^{R}_{K}(\Theta(x_1))x_0).
\end{equation}
Therefore $\Theta$ is a solution of the strong Mauer-Cartan equation
by Corollary 5.2.
\end{proof}

In view of Proposition 4.4, Theorem 5.4 and Proposition 4.7, the
following result is clear.

 {\bf Theorem 5.6.} Suppose that $K$ is
 an
$\mathcal{O}$-operator on a Leibniz algebra $(\mathfrak g, [ \ , \
])$ with respect to a representation $(V, \rho^{L},\rho^{R})$
 and
$\Theta:\mathfrak g \longrightarrow V$ is a solution of the strong
Maurer-Cartan equation on the twilled Leibniz algebra $\mathfrak
g\bowtie V_{K}$. Then $K\Theta K$ is an $\mathcal{O}$-operator on
$(\mathfrak g, [ \ , \ ])$ with respect to $(V, \rho^{L},\rho^{R})$.
Moreover, $K$ and $K\Theta K$ are compatible.

\section{Leibniz
$r-n$ structures, RBN-structures and $\mathcal{B}N$-structures }

\quad {\bf Definition 6.1.} Let $\pi$ be a classical Leibniz
$r$-matrix and $N: \mathfrak g\longrightarrow \mathfrak g$ a
Nijenhuis operator on a Leibniz algebra $(\mathfrak g,[ \ , \ ])$.
We say that $(\pi,N)$ is a {\it Leibniz $r-n$ structure} on
$(\mathfrak g,[ \ , \ ])$ if
for any 
$\alpha,\beta\in \mathfrak g^{*}$, 
\begin{align}\label{7.1}
N\pi^{\sharp}&=\pi^{\sharp}N^{*},\\ \label{7.2}
[\alpha,\beta]^{N\pi^{\sharp}}&=[\alpha,
\beta]_{N^{*}}^{\pi^{\sharp}},
\end{align}
where
$[\alpha,\beta]^{N\pi^{\sharp}}=L^{*}(N\pi^{\sharp}(\alpha))(\beta)-((L^{*}+R^{*})(N\pi^{\sharp}(\beta)))(\alpha)$.

 We immediately have the following.

{\bf Theorem 6.2.} Suppose
$(\pi,N)$ is a Leibniz $r-n$ structure on a Leibniz algebra
$(\mathfrak g,[ \ , \ ])$, then $(\pi^{\sharp},  N^{*}, N)$ is a
dual $\mathcal{O}$N-structure on $(\mathfrak g,[ \ , \ ])$
associated to the dual representation $(\mathfrak
g^{*},L^{*},-L^{*}-R^{*})$.

{\bf Definition 6.3.} Let $\mathcal{R}:\mathfrak g\longrightarrow
\mathfrak g$ be a Rota-Baxter operator and $N: \mathfrak
g\longrightarrow \mathfrak g$ a Nijenhuis operator on a Leibniz
algebra $(\mathfrak g,[ \ , \ ])$. A pair $(\mathcal{R},N)$ is
called an {\it RBN-structure} on $(\mathfrak g,[ \ , \ ])$ if
\begin{align}\label{7.3}
N\mathcal{R}&=\mathcal{R}N,\\ \label{7.4} [x,
y]^{N\mathcal{R}}&=[x,y]_{N}^{\mathcal{R}},~~~ \forall ~x,~y\in
\mathfrak g,
\end{align}
where $[x,
y]^{N\mathcal{R}}=[N\mathcal{R}(x),y]+[x,N\mathcal{R}(y)].$

 {\bf Example 6.4.} Let $(\mathfrak g,[ \ , \ ])$ be the Leibniz
algebra with basis $\{\varepsilon_{1},\varepsilon_{2}\}$ and
 the nonzero multiplication be
 $$[\varepsilon_{1}, \varepsilon_{1}]=\varepsilon_{2}.$$
Define linear maps $\mathcal{R},N: \mathfrak g\longrightarrow \mathfrak g$
respectively by the matrices
$$\mathcal{R}(\varepsilon_{1}, \varepsilon_{2})= (\varepsilon_{1}, \varepsilon_{2})\left(
  \begin{array}{ccc}
    0& 0 \\
    a_{21} & a_{22} \\
  \end{array}
\right)$$ and $$ N(\varepsilon_{1}, \varepsilon_{2})=
(\varepsilon_{1}, \varepsilon_{2})\left(
  \begin{array}{ccc}
    a& 0 \\
    0 & a \\
  \end{array}
\right)
$$  with respect to the basis $\{\varepsilon_1, \varepsilon_2
\}$. Then $\mathcal{R}$ is a Rota-Baxter operator and $N: \mathfrak
g\longrightarrow \mathfrak g$ is a Nijenhuis operator on $\mathfrak
g$, and direct computation shows that $(\mathcal{R},N)$ is an RBN-structure.


 {\bf Definition 6.5.} \cite{ST} A {\it quadratic Leibniz algebra} is a Leibniz algebra $(\mathfrak g,[ \ , \ ])$ with a
nondegenerate skew-symmetric bilinear form $\mathfrak{q}\in
\wedge^{2}\mathfrak g^{*}$ satisfying the following invariant
condition:
$$\mathfrak{q}(x_0, [x_1, x_2]) =  \mathfrak{q}([x_0, x_2] + [x_2, x_0],
x_1),~~ \forall~ x_0, x_1, x_2 \in \mathfrak g.$$

Let $(\mathfrak g,[ \ , \ ],\mathfrak{q})$ be a quadratic Leibniz
algebra. Then
$$({\mathfrak{q}}^{\sharp})^{-1}:\mathfrak
g\longrightarrow \mathfrak g^{*},~\langle
({\mathfrak{q}}^{\sharp})^{-1}(x),y\rangle=\mathfrak{q}(x,y),~~\forall~
x,y\in \mathfrak g$$
 is an isomorphism from the regular
representation $(\mathfrak g, L, R)$ to its dual representation
$(\mathfrak g^{*},L^{*},-R^{*}-L^{*})$. Moreover,
\begin{equation}\label{7.5}
({\mathfrak{q}}^{\sharp})^{-1}R(x)=(-R^{*}-L^{*})(x)({\mathfrak{q}}^{\sharp})^{-1},
~~~({\mathfrak{q}}^{\sharp})^{-1}L(x)=L^{*}(x)({\mathfrak{q}}^{\sharp})^{-1},~~\forall~
x\in \mathfrak g.
\end{equation}

 {\bf Theorem 6.6.} Let $(\mathfrak g,[ \ , \ ],\mathfrak{q})$ be a quadratic Leibniz
 algebra, $\mathcal{R}:\mathfrak g\longrightarrow \mathfrak g$ a linear map, and
  $N:\mathfrak g\longrightarrow \mathfrak g$ a Nijenhuis operator on $(\mathfrak g,[ \ , \ ])$.
  Suppose that $\pi^{\sharp}=\mathcal{R}\mathfrak{q}^{\sharp}$ and $\pi$ is symmetric such that
\begin{equation}\label{7.6}
\mathfrak{q}^{\sharp}N^{*}=N\mathfrak{q}^{\sharp}.
\end{equation}
Then $(\mathcal{R},N)$ is an RBN-structure on $(\mathfrak g,[\ ,\
],\mathfrak{q})$ if and only if $(\pi,N)$ is a Leibniz $r-n$
structure.

\begin{proof} Since
 $\mathfrak{q}^{\sharp}$ is bijective,
 for any $\alpha_0,\alpha_1\in
 \mathfrak g^{*}$ there are $x_0,x_1\in \mathfrak g$ such that
 $$\alpha_0=(\mathfrak{q}^{\sharp})^{-1}(x_0),\alpha_1=(\mathfrak{q}^{\sharp})^{-1}(x_1).$$
By \cite[Corollary 4.17]{ST}, $\pi$ is a Leibniz
 $r$-matrix if and only if $\pi^{\sharp}(\mathfrak{q}^{\sharp})^{-1}$ is a Rota-Baxter
 operator on $(\mathfrak g,[ \ , \ ])$.

$(\Longrightarrow)$ By \eqref{7.3} and \eqref{7.6}, we have
\begin{equation*}
N\pi^{\sharp}=\pi^{\sharp}N^{*}.
\end{equation*}
We now verify that
$$[\alpha_0, \alpha_1]^{N\pi^{\sharp}}
=[\alpha_0,\alpha_1] _{N^{*}}^{\pi^{\sharp}}.$$
In fact, by \eqref{7.5} we have
\begin{eqnarray}\nonumber
[\alpha_0, \alpha_1]^{\pi^{\sharp}}
&=&L^{*}(\pi^{\sharp}(\alpha_0))\alpha_1-(L^{*}+R^{*})(\pi^{\sharp}(\alpha_1))\alpha_0\\
\nonumber
&=&L^{*}(\pi^{\sharp}(\mathfrak{q}^{\sharp})^{-1}(x_0))(\mathfrak{q}^{\sharp})^{-1}(x_1)
-(L^{*}+R^{*})(\pi^{\sharp}(\mathfrak{q}^{\sharp})^{-1}(x_1))(\mathfrak{q}^{\sharp})^{-1}(x_0)\\
\nonumber
&=&L^{*}(\mathcal{R}(x_0))(\mathfrak{q}^{\sharp})^{-1}(x_1)-(L^{*}+R^{*})(\mathcal{R}(x_1))(\mathfrak{q}^{\sharp})^{-1}(x_0)\\
\nonumber
&=&(\mathfrak{q}^{\sharp})^{-1}(L(\mathcal{R}(x_0))(x_1)+R(\mathcal{R}(x_1))(x_0))\\
\label{7.7}
&=&(\mathfrak{q}^{\sharp})^{-1}([\mathcal{R}(x_0) ,x_1]+[x_0, \mathcal{R}(x_1)])=(\mathfrak{q}^{\sharp})^{-1}[x_0,x_1]^{\mathcal{R}}.   
\end{eqnarray}
According to \eqref{7.6} and \eqref{7.7},
\begin{eqnarray*}
&&[\alpha_0, \alpha_1]^{N\pi^{\sharp}} -[\alpha_0,
\alpha_1]_{N^{*}}^{\pi^{\sharp}}=[\alpha_0,
\alpha_1]^{N\pi^{\sharp}} -[N^{*}(\alpha_0),
\alpha_1]^{\pi^{\sharp}}-[\alpha_0,
N^{*}(\alpha_1)]^{\pi^{\sharp}}+N^{*}[\alpha_0,\alpha_1] ^{\pi^{\sharp}}\\
&=&(\mathfrak{q}^{\sharp})^{-1}([x_0,x_1]^{N\mathcal{R}}
)-[N^{*}(\mathfrak{q}^{\sharp})^{-1}(x_0),(\mathfrak{q}^{\sharp})^{-1}(x_1)]^{\pi^{\sharp}}\\&&-
  [(\mathfrak{q}^{\sharp})^{-1}(x_0),N^{*}(\mathfrak{q}^{\sharp})^{-1}(x_1)]^{\pi^{\sharp}} +
N^{*}[(\mathfrak{q}^{\sharp})^{-1}(x_0),(\mathfrak{q}^{\sharp})^{-1}(x_1)]^{\pi^{\sharp}}
\\
&=&(\mathfrak{q}^{\sharp})^{-1}[x_0,x_1]^{N\mathcal{R}}
-[(\mathfrak{q}^{\sharp})^{-1}(N(x_0)),(\mathfrak{q}^{\sharp})^{-1}(x_1)]^{\pi^{\sharp}}
-[(\mathfrak{q}^{\sharp})^{-1}(x_0),(\mathfrak{q}^{\sharp})^{-1}N(x_1)]^{\pi^{\sharp}}
\\&&+ N^{*} [(\mathfrak{q}^{\sharp})^{-1}(x_0),(\mathfrak{q}^{\sharp})^{-1}(x_1)]^{\pi^{\sharp}}
\\
&=&(\mathfrak{q}^{\sharp})^{-1}([x_0,x_1]^{N\mathcal{R}}
)-(\mathfrak{q}^{\sharp})^{-1}([N(x_0),x_1]^{\mathcal{R}}
)-(\mathfrak{q}^{\sharp})^{-1}[x_0, N(x_1)]^{\mathcal{R}}+ N^{*}
(\mathfrak{q}^{\sharp})^{-1}([x_0,x_1]^{\mathcal{R}})
\\
&=&(\mathfrak{q}^{\sharp})^{-1}[x_0,x_1]^{N\mathcal{R}}
-(\mathfrak{q}^{\sharp})^{-1}[N(x_0),x_1]^{\mathcal{R}}
-(\mathfrak{q}^{\sharp})^{-1}[x_0,N(x_1)]^{\mathcal{R}}+
(\mathfrak{q}^{\sharp})^{-1}N [x_0,x_1]^{\mathcal{R}}
\\
&=&(\mathfrak{q}^{\sharp})^{-1}[x_0,x_1]^{N\mathcal{R}}
-(\mathfrak{q}^{\sharp})^{-1}([N(x_0),x_1]^{\mathcal{R}} +[x_0,
N(x_1)]^{\mathcal{R}} - N([x_0,x_1]^{\mathcal{R}})
\\
&=&(\mathfrak{q}^{\sharp})^{-1}([x_0,x_1]^{N\mathcal{R}}
-[x_0,x_1]_{N}^{\mathcal{R}})=0.
\end{eqnarray*}
Therefore
$$[\alpha_0,\alpha_1]^{N\pi^{\sharp}}  =[\alpha_0,\alpha_1]
_{N^{*}}^{\pi^{\sharp}}.$$ Hence $(\pi,N)$ is a Leibniz $r-n$
structure.

 $(\Longleftarrow)$ By \eqref{7.1} and \eqref{7.6}, we get
 $$N\mathcal{R}=\mathcal{R}N.$$
The remaining part is similar to the converse argument.
\end{proof}

A symmetric nondegenerate bilinear form $\mathcal{B}\in
\wedge^{2}\mathfrak g^{*} $ on a Leibniz algebra $(\mathfrak g,[ \ ,
\ ])$ induces a linear map $\mathcal{B}^{\sharp}: \mathfrak
g^{*}\longrightarrow \mathfrak g$ by
$$\langle (\mathcal{B}^{\sharp})^{-1}(x_0),x_1\rangle=\mathcal{B}(x_0,x_1),~~\forall ~x_0,~x_1\in \mathfrak g.$$
Here $\mathcal{B}$ is nondegenerate if $\mathcal{B}^{\sharp}:
\mathfrak g^{*}\longrightarrow \mathfrak g$ is an isomorphism.

 {\bf Definition 6.7.} Let $N$ be a Nijenhuis operator on a Leibniz algebra $(\mathfrak g,[ \ , \ ])$ and
$\mathcal{B}\in \wedge^{2}\mathfrak g^{*}$ a symmetric nondegenerate
bilinear form satisfying the closeness condition:
\begin{equation}\label{7.8}
\mathcal{B}(x_2,[x_0,x_1])=-\mathcal{B}(x_1,[x_0,x_2])+\mathcal{B}(x_0,[x_1,x_2])+\mathcal{B}(x_0,[x_2,x_1])
\end{equation}
for any $x_0,x_1,x_2\in \mathfrak g$.
 We say that $(\mathcal{B}, N)$ is
a {\it $\mathcal{B}N$-structure} on $\mathfrak g$ if
\begin{equation}\label{7.9}\mathcal{B}(N(x_0),x_1)=\mathcal{B}(x_0,N(x_1))
\end{equation}
and
\begin{equation}\label{7.10}
\mathcal{B}(x_2,N[x_0,x_1])=-\mathcal{B}(x_1,N[x_0,x_2])+\mathcal{B}(x_0,N[x_1,x_2])+\mathcal{B}(x_0,N[x_2,x_1]).
\end{equation}

 {\bf Theorem 6.8.} Let $N$ be a
Nijenhuis operator on a Leibniz algebra$(\mathfrak g,[ \ , \ ])$ and
$\mathcal{B}\in \wedge^{2}\mathfrak g^{*}$
 a nondegenerate symmetric bilinear form satisfying identity (6.8). If
$(\mathcal{B}, N)$ is a $\mathcal{B}N$-structure, then
 $(\mathcal{B}^{\sharp}, S=N^{*},N)$ is a dual
$\mathcal{O}$N-structure on $\mathfrak g$ associated to the dual
representation $(\mathfrak g^{*},L^{*}, -R^{*}-L^{*})$.
\begin{proof} In view of (6.9),
 $\mathcal{B}^{\sharp}N^{*}=N\mathcal{B}^{\sharp}.$ Since
 $\mathcal{B}$ is symmetric and satisfies (6.8),
 $\mathcal{B}^{\sharp} $ is an $\mathcal{O}$-operator on $(\mathfrak g,[ \ , \ ])$ with
  respect to the representation $(\mathfrak g^{*},L^{*}, -R^{*}-L^{*})$.
  Now we claim that 
\begin{equation}\label{7.11}
[\alpha_0,\alpha_1]^{N\mathcal{B}^{\sharp}}=[\alpha_0,\alpha_1]_{N^{*}}^{\mathcal{B}^{\sharp}}.
\end{equation}
Since
 $\mathcal{B}^{\sharp}$ is bijective, write
 $\alpha_0=(\mathcal{B}^{\sharp})^{-1}(x_0),\alpha_1=(\mathcal{B}^{\sharp})^{-1}(x_1)$ for some $x_0,x_1\in \mathfrak g$.
For any $x_2\in \mathfrak g$, based on \eqref{7.8} and \eqref{7.9},
\begin{eqnarray*}&&\langle[\alpha_0,\alpha_1]^{N\mathcal{B}^{\sharp}}
-[\alpha_0,\alpha_1]_{N^{*}}^{\mathcal{B}^{\sharp}},x_2\rangle
\\&=&\langle
L^{*}(N\mathcal{B}^{\sharp}(\alpha_0))(\alpha_1)-((L^{*}+R^{*})(N\mathcal{B}^{\sharp}(\alpha_1)))(\alpha_0)
-L^{*}(\mathcal{B}^{\sharp}N^{*}(\alpha_0)(\alpha_1)
\\&&+((L^{*}+R^{*})(\mathcal{B}^{\sharp}(\alpha_1)))(N^{*}(\alpha_0))
-L^{*}(\mathcal{B}^{\sharp}(\alpha_0)(N^{*}(\alpha_1))
\\&&+((L^{*}+R^{*})(\mathcal{B}^{\sharp}N^{*}(\alpha_1)))(\alpha_0)
+N^{*}L^{*}(\mathcal{B}^{\sharp}(\alpha_0)(\alpha_1)
-N^{*}((L^{*}+R^{*})(\mathcal{B}^{\sharp}(\alpha_1)))(\alpha_0),x_2\rangle
\\&=&-\langle (\mathcal{B}^{\sharp})^{-1}(x_1),[N(x_0),x_2]\rangle
+\langle
(\mathcal{B}^{\sharp})^{-1}(x_0),[N(x_1),x_2]+[x_2,N(x_1)]\rangle
\\&&+\langle (\mathcal{B}^{\sharp})^{-1}(x_1),[N(x_0),x_2]\rangle
-\langle
(\mathcal{B}^{\sharp})^{-1}(N(x_0)),[x_1,x_2]+[x_2,x_1]\rangle
\\&&
+\langle (\mathcal{B}^{\sharp})^{-1}(N(x_1)),[x_0,x_2]\rangle
-\langle
(\mathcal{B}^{\sharp})^{-1}(x_0),[N(x_1),x_2]+[x_2,N(x_1)]\rangle
\\&&-\langle (\mathcal{B}^{\sharp})^{-1}(x_1),[x_0,N(x_2)]\rangle
+\langle
(\mathcal{B}^{\sharp})^{-1}(x_0),[x_1,N(x_2)]+[N(x_2),x_1]\rangle
\\&=&-\mathcal{B}(x_1,[N(x_0),x_2])
+\mathcal{B}(x_0,[N(x_1),x_2]+[x_2,N(x_1)])
+\mathcal{B}(x_1,[N(x_0),x_2])
\\&&-\mathcal{B}(N(x_0),[x_1,x_2]+[x_2,x_1])
+\mathcal{B}(N(x_1),[x_0,x_2])
-\mathcal{B}(x_0,[N(x_1),x_2]+[x_2,N(x_1)])
\\&&-\mathcal{B}(x_1,[x_0,N(x_2)])
+\mathcal{B}(x_0,[x_1,N(x_2)]+[N(x_2),x_1])
\\&=&
-\mathcal{B}(N(x_0),[x_1,x_2]+[x_2,x_1])+\mathcal{B}(N(x_1),[x_0,x_2])
\\&&-\mathcal{B}(x_1,[x_0,N(x_2)])
+\mathcal{B}(x_0,[x_1,N(x_2)]+[N(x_2),x_1])
\\&=&
-\mathcal{B}(x_0,N([x_1,x_2]+[x_2,x_1]))+\mathcal{B}(x_1,N([x_0,x_2]))
+\mathcal{B}(N(x_2),[x_0,x_1])
 \\&=&0,
\end{eqnarray*}
which implies that \eqref{7.11}  holds.
\end{proof}

The following result is a consequence of Proposition 4.7 and Theorem
6.8.

{\bf Corollary 6.9.} Let $N$ be a Nijenhuis operator on a Leibniz
algebra $(\mathfrak g, [ \ , \ ])$ and $\mathcal{B}\in
\wedge^{2}\mathfrak g^{*}$ a nondegenerate symmetric bilinear form.
If $(\mathcal{B}, N)$ is a $\mathcal{B}N$-structure, then
$\mathcal{B}^{\sharp}$ and $\mathcal{B}^{\sharp}N^{*}$ are
compatible $\mathcal{O}$-operators on $(\mathfrak g, [ \ , \ ])$
with respect to $(\mathfrak g^{*},L^{*}, -R^{*}-L^{*})$.

\section*{Acknowledgments}
\quad The work is supported by the National Natural Science
Foundation of China (grant no. 12171303), Simons Foundation (grant
no. 523868), the Natural Science Foundation of Zhejiang Province of
China (grant no. LY19A010001) and the Science and Technology
Planning Project of Zhejiang Province (2022C01118).

\end{document}